\documentclass[final,onefignum,onetabnum]{siamart190516}

\usepackage{lipsum}
\usepackage{amsfonts}
\usepackage{amsmath, amssymb, textcomp}
\usepackage{graphicx}
\usepackage{subcaption}
\usepackage{epstopdf}
\usepackage{algorithmic}
\usepackage{bm}
\usepackage{color}

\ifpdf
  \DeclareGraphicsExtensions{.eps,.pdf,.png,.jpg}
\else
  \DeclareGraphicsExtensions{.eps}
\fi

\usepackage[section]{placeins}

\newsiamremark{remark}{Remark}
\newsiamremark{hypothesis}{Hypothesis}
\crefname{hypothesis}{Hypothesis}{Hypotheses}
\newsiamthm{claim}{Claim}
\graphicspath{{./Graphs/}}

\renewcommand\thesection{\arabic {section}}

\newtheorem{defn}{Definition}[section]

\theoremstyle{plain} 

\theoremstyle{plain} 

\theoremstyle{plain}

\hypersetup{urlcolor=blue, citecolor=red}

\newcommand{\ind}{\operatorname{ind}}
\renewcommand{\d}{\mathrm{d}}
\newcommand{\ir}{\operatorname{ind}_\textrm{rad}}
\newcommand{\ia}{\operatorname{ind}_\textrm{ang}}
\newcommand{\floor}{\operatorname{floor}}
\newcommand{\ceil}{\operatorname{ceil}}
\DeclareMathOperator{\MET}{MET}
\DeclareMathOperator{\MIX}{MIX}

\headers{Modelling freezing of gait}
{A Wang, J Sieber, W Young and K Tsaneva-Atanasova } 

\title{Time series analysis and modelling of the freezing of gait phenomenon}

\author{Ai Wang\thanks{Department of Mathematics and Statistics,  University of Exeter, North Park Road, Exeter EX4 4QL, UK}
\and Jan Sieber\footnotemark[1]
 \and William R. Young\thanks{Sport \& Health Sciences, St Luke's Campus, Heavitree Road, Exeter EX1 2LU UK}
\and Krasimira Tsaneva-Atanasova\footnotemark[1] \thanks{Living Systems Institute, EPSRC Hub for Quantitative Modelling in Healthcare, University of Exeter, Stocker Road, Exeter EX4 4QD, UK and Institute for Advanced Study, Technical University of Munich, Lichtenbergstrasse 2 a, D-85748 Garching, Germany.  (\email{K.Tsaneva-Atanasova@exeter.ac.uk}.)}}

\usepackage{amsopn}

\begin{document}

\maketitle
\begin{abstract}
  Freezing of Gait (FOG) is one of the most debilitating symptoms of Parkinson's Disease and is associated with falls and loss of independence. The patho-physiological mechanisms underpinning FOG are currently poorly understood. In this paper we combine time series analysis and mathematical modelling to study the FOG phenomenon's dynamics. We focus on the transition from stepping in place into freezing and treat this phenomenon in the context of an escape from an oscillatory attractor into an equilibrium attractor state. We extract a discrete-time discrete-space Markov chain from experimental data and divide its state space into communicating classes to identify the transition into freezing. This allows us to develop a methodology for computationally estimating the time to freezing as well as the phase along the oscillatory (stepping) cycle of a patient experiencing Freezing Episodes (FE). The developed methodology is general and could be applied to any time series featuring transitions between different dynamic regimes including time series data from forward walking in people with FOG .
\end{abstract}



\section{Introduction}
People with Parkinson's Disease (PD) will often walk with reduced gait speed and shorter stride length. Spatial and temporal characteristics of Parkinsonian gait are also typically more variable compared to age-matched controls~\cite{PMID:9613733,hausdorff2005gait}. Furthermore, approximately 50\% of people with advanced Parkinson's will experience freezing of gait (FOG)~\cite{PMID:11261746}. Patients describe FOG as the sensation that their feet are glued to the floor, preventing them from initiating a new step. Indeed, a greater variability in walking patterns has been observed in Parkinson's patients with FOG. This is characterised by increased step coefficient of variation~\cite{Almeida513}, asymmetry, rhythmicity~\cite{Nantel2011} and difficulty  coordinating~\cite{Plotnik11}, compared to patients without freezing. These studies focus on descriptive statistics of the walking and/or stepping time series. In this study, in addition to performing advanced time series analysis we concentrate on revealing the dynamic (geometric) properties of the transitions between stepping and freezing.

Dynamic modelling and analysis can help in understanding and characterising specific features and properties of Parkinsonian gait. This in turn could inform future rehabilitation and prevention interventions as well as strategies that people with Parkinson's might benefit from through informal use in daily life. A variety of mathematical modelling and data analysis approaches have been applied in the context of Parkinson's gait and motor control more generally as reviewed in \cite{sarbaz2016review}. 
Data-driven prediction and detection of the FOG phenomenon have been extensively addressed in previous work, including by us in the context of stepping in place force platform data \cite{parakkal2020data},  and others (see Table 1 in \cite{pardoel2019wearable} for list of studies based on wearable sensors data \cite{pardoel2019wearable}), and more recently plantar pressure data \cite{pardoel2021grouping,shalin2021prediction}. A stochastic model of gait consisting of a random walk on a chain has been also proposed and applied to describe alterations in gait dynamics from childhood to adulthood \cite{ashkenazy2002stochastic}. 
To our knowledge none of the studies to date have combined mathematical modelling and data analysis to investigate dynamic properties of the FOG phenomenon in Parkinson's Disease, which is the main focus of our work.

\paragraph{Summary of results}
We analyse the transition between stepping and freezing observed during stepping-in-place experiments performed by  Nantel et al \cite{Nantel2011}. In our analysis, after phase space embedding, stepping motion appears as large-amplitude oscillations while freezing appears as an approximate equilibrium or irregular drifting. 
We explore whether there exists a  preferred phase at which the stepping trajectory escapes from the regular periodic stepping behavior to transition into freezing. The first step toward an answer to this question is to identify for each freezing event a location in phase space (after embedding) at which the transition occurs. At first sight this appears to depend strongly on the parameter choices in our data processing, such as, for example, a threshold for freezing that we could define. To address this sensitivity we seek to develop a recurrence-based method for identifying a unique phase and amplitude of the point at which the transition occurs that does not depend on method parameter choices. Our method consists of three steps. First, we apply a Hilbert Transform to reconstruct a two-dimensional embedding from the scalar experimental times series data as a signal in the complex plane. Second, we construct a Markov chain by discretizing the complex plane and counting the empirical transition probabilities of the Hilbert Transform output between the boxes from the discretization. Over a large range of method parameters the Markov chain has a clearly identifiable largest communicating class corresponding to regular periodic stepping (the \emph{stepping class}.) In the final step we determine which boxes have a large transition probability out of the stepping class. We take the location (phase and amplitude) of these boxes in the complex plane as the threshold for transition to freezing. With the help of this naturally emerging threshold we obtain a mean escape time characterising the transition from stepping to freezing. We can also investigate if these transitions occur independent of phase or in some range of preferred phases (angles) along the oscillatory cycle preceding the freezing episodes.
	
The paper presents the method for identifying the transition's phase and amplitude for each event, demonstrates that it is robust with respect to choice of method parameters, and tests it on the small patient data set that we have currently available. In principle, pinpointing the transition will allow us to determine whether transitions occur preferably at universal or patient-specific phase ranges, or whether they are uniformly distributed (we treat the latter as the null hypothesis). This would help us to further understand the underlying dynamics of FOG.

The paper is organized as follows. After describing the type of data we use for our method development in \Cref{sec:data:intro},  we present our null model of the transition between stepping and freezing in \Cref{sec:nullmodel}. \Cref{sec:method,sec:mc:prop} then describe our method to construct a Markov chain from the data and the analysis of the Markov chain's properties we expect for our data. The new algorithm we develop is applied to patients' stepping data in \Cref{sec:data} and our findings regarding phase dependence of the transition into freezing are given in \Cref{sec:phase:prediction}. \Cref{sec:conclusions} follows with conclusions and outlook.

\section{Stepping-in-place experiments}
\label{sec:data:intro}
Nantel \emph{et al}.\ \cite{Nantel2011} performed experiments
recording freezing of gait using a repetitive stepping-in-place task
on force plates which could identify freezing episodes in subjects
with Parkinson's Disease.  The task consisted of alternatingly raising
the legs at a self-selected pace for $90$\,s per trial. Ground
reaction forces were measured with a sampling frequency of $100$\,Hz
on two force plates. The series of experiments generated 6 trials'
data sets of $90$\,s force time series for each of the 34 subjects.

\begin{figure}[ht]
	\centering
	\setlength{\abovecaptionskip}{0.cm}
	\setlength{\belowcaptionskip}{-0.cm}
	\includegraphics[width=\textwidth]{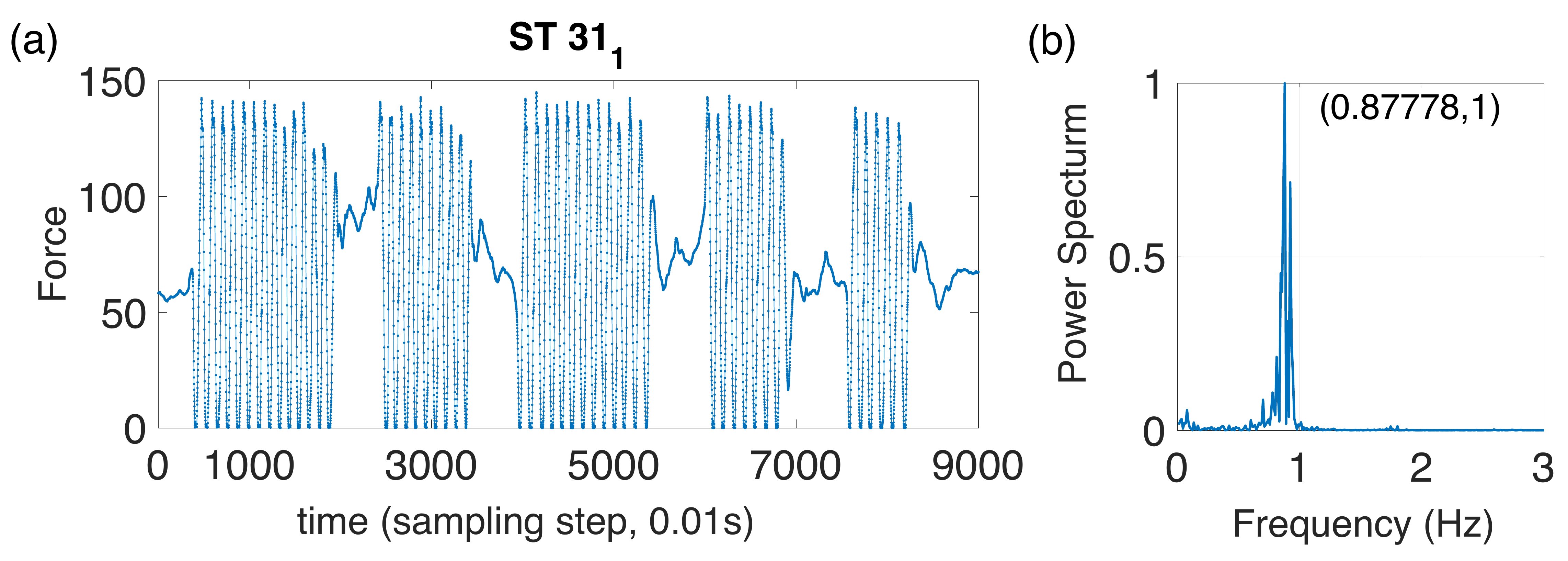}
	\caption{Panel (a): Time series of the left foot vertical force (as a percentage of body weight) of the subject numbered ST 31, data set $1$, sampling time step $\delta t=0.01$; Panel (b): Power spectrum (scaled
                  to maximum equal to unit) with strong dominant peak
                  at $f_\mathrm{stp}\approx 0.88$Hz;}
	\label{p1}
\end{figure}
Figure \ref{p1} shows a typical time series generated by the data
collection. The time series represent the left foot vertical force, given as a percentage of the body weight,
 of one of the subjects numbered ST\,31. The
sampling time step $\delta t$ equals 0.01 seconds. We observe several
sudden transitions in the force magnitude from a large-amplitude
oscillatory behavior to small-amplitude irregular fluctuations and
drifting at approximately half the vertical force, which indicates
that both feet are on the ground and the patient
experiences a freezing episode. Episodes of standstill (involuntary freezing or intentional) are uniquely identified by absence of zero-force periods for longer than the identified stepping frequency shown in \cref{p1}(b). An expert manually identified
freezing episodes in the time series data during the experiment in order to distinguish freezing episodes from episodes when the subject is
intentionally standing still, as seen at the start and end of many
time series.

During our analysis we non-dimensionalize the force data to the interval $[-1,1]$,  shifting the freezing to approximately the value $0$, by applying the transformation
\begin{align}\label{force_scaling}
  \mathrm{Force}_\mathrm{scal}(t)=\frac{\mathrm{Force}(t)-\operatorname{mean}(\mathrm{Force}(t))}{\max|\mathrm{Force}(t)-\operatorname{mean}(\mathrm{Force}(t))|}\mbox{.}
\end{align}

We report time in multiples of the sampling time step $\delta t$, such
that all times $t$ are integers. When presenting experimental data in
the paper we label subject number (ST), data set and time interval (in
multiples of $\delta t=0.01$). For example, \Cref{p1} shows data from
ST\,$31$, data set $1$, and the entire time interval $[0,9000]$
(corresponding to the full $90$\,s).  In the force measurements,
stepping is seen as a regular periodic behavior similar to a limit
cycle oscillation, while freezing resembles an equilibrium or
irregular drifting.

Stepping
in place has been used as a protocol because it is extremely
challenging to evaluate FOG using traditional force platforms, due to
their limited size and the practical difficulties of observing
unpredictable FOG events when asking participants to walk over
them. For this reason, the approach of \cite{Nantel2011} has been to ask
participants to step in place when standing on two force
platforms. This approach provides highly accurate dynamic force
information for each foot separately. Because these data relate to a protocol where patients have been denied the opportunity to walk forward, the temporal aspects of the task are more prominent (spatial features such as stride length are absent). Therefore, we focus our analysis on temporal features, i.e dynamical properties of scalar force data.

Full data sets and processing scripts are
available at the following link \url{https://figshare.com/s/a14be7360925639736ba}.

\section{Nullmodel: transitions independent of phase}\label{sec:nullmodel}
\begin{figure}[ht]
		\centering
 		\setlength{\abovecaptionskip}{0.cm}
 		\setlength{\belowcaptionskip}{-0.cm}
 		\includegraphics[width=0.9\textwidth]{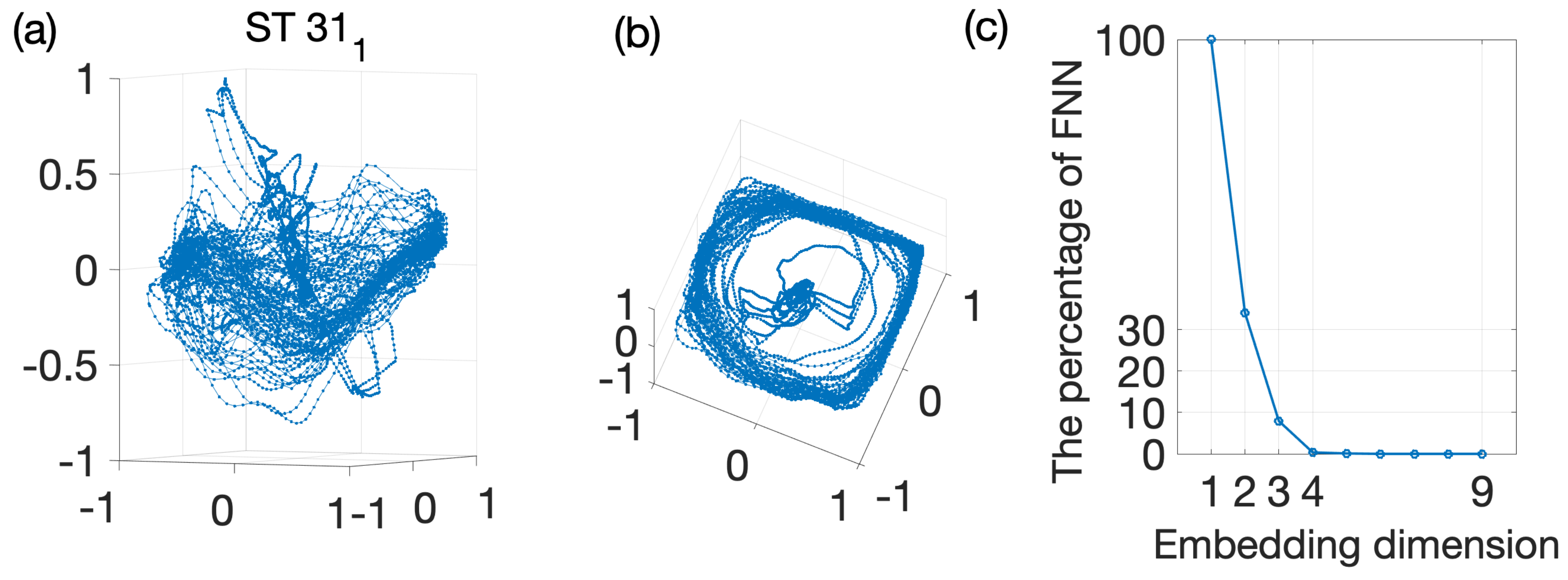}
 		\caption{ Time delay
                  embedding analysis of time profile shown in
                  \cref{p1}: Panel (a) and (b): Delay embedding in 3 dimensional space; Panel (c): False nearest neighbor percentage for different
                  embedding dimensions and time delay.}
 		\label{delembedding}
 	\end{figure}

Time profiles during regular stepping exhibit strong periodicity, as
\cref{p1}(a) shows. \Cref{p1}(b) 
depicts the corresponding power spectrum characterized by a distinct dominant peak,
corresponding to the stepping frequency $f_\mathrm{stp}$ of the
subject during the experiment (in this case
$f_\mathrm{stp}\approx 0.88$Hz, such that the angular stepping
frequency is $\omega_\mathrm{stp}\approx1.76\pi/s$). The False Nearest
Neighbor (FNN) criterion \cite{kennel1992determining} for time
delay embedding with $\Delta_t$ equal to a quarter of the
dominant period, $T_\mathrm{stp}/4=\pi/(2\omega_\mathrm{stp})$, demonstrates
(as seen in \cref{delembedding}(right)) that two dimensions predict 65\% of
the signal correctly (that is, $35\%$ of the predictions from nearest
neighbours mismatch according to the FNN criterion), while
10\% is mismatched with three embedding dimensions. The third
embedding dimension plays a role after freezing events have occurred,
as the time-delay embedded phase portraits in panel (a) and (b) of
\cref{delembedding} illustrate.

\begin{figure}[h]
		\centering
 		\setlength{\abovecaptionskip}{0.cm}
 		\setlength{\belowcaptionskip}{-0.cm}
 		\includegraphics[width=\textwidth]{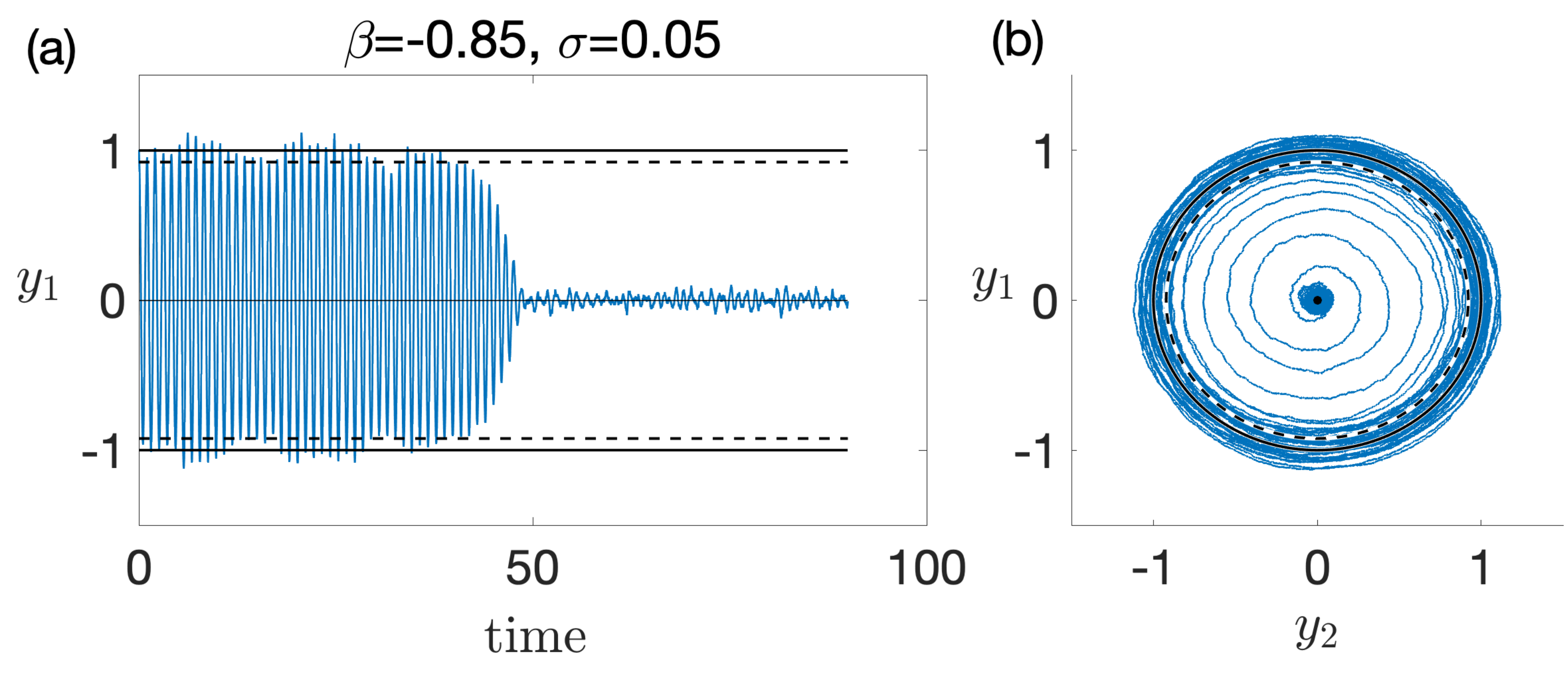}
 		\caption{ An example realization of the noise-driven dynamics of (\eqref{eq:HNF}). Both Panel show the same realization in blue for $\beta = -0.85$ and $\sigma = 0.05$  in the phase space, for which the deterministic part contains a stable equilibrium at the origin, an unstable limit cycle (black dotted line), and stable limit cycle (solid line). }
 		\label{modelsimulation}
 	\end{figure}

The escaping-from-a-limit-cycle behavior visible in \cref{delembedding} could be accounted for (phenomenologically) by bi-stable oscillatory dynamics featuring co-existence of an attracting periodic limit cycle, corresponding to periodic stepping, and an equilibrium state, corresponding to freezing. Noisy fluctuations and/or external inputs could drive transitions between these two stable states representing transitions from stepping into freezing.
The generalized Hopf normal form model \cite{Guckenheimer:2007} exhibits such dynamics and hence can be used for describing some aspects of the behavior shown by the time-delay embedding of the data \cref{delembedding}. Accordingly, the model we consider is based on the generalized Hopf normal form model with additive white noise and has the following form:
\begin{subequations}
		\begin{alignat}{2}
		\d y_1(t) &=\left[ \beta y_1-2\pi\omega y_2+(1-\beta)y_1(y_1^2+y_2^2)-y_1(y_1^2+y_2^2)^2\right] \d t + \sigma \d W_1 (t), \\
		\d y_2(t) &= \left[ 2\pi\omega y_1+\beta y_2+(1-\beta)y_2(y_1^2+y_2^2)-y_2(y_1^2+y_2^2)^2\right] \d t + \sigma \d W_2 (t),
		\end{alignat}
	\label{eq:HNFc}
\end{subequations}
where  the initial condition is $y_1= 1,y_2 = 0$, and $W_j(t)$ are standard Wiener processes. Here $\sigma$ is the noise amplitude and $\sigma^2$ is the variance with $\sigma> 0$. 
The stochastic differential equation \eqref{eq:HNFc} can be transformed into polar coordinates by $z(t) = R(t)\exp[i\theta(t)]$, resulting in  \cite{Creaser2018}:
	\begin{subequations}
		\begin{alignat}{2}
		\d R&= \left[ \beta R+(1-\beta)R^3-R^5+\dfrac{\sigma^2}{2R}\right]\d t+\sigma \d W_R, \label{rpolar}\\
		\d\theta & = 2\pi\omega \d t+ \dfrac{\sigma}{R}\d W_\theta, 	\label{thetapolar}
		\end{alignat}
			\label{eq:HNF}
	\end{subequations}
where $W_R$ and $W_\theta$ are independent standard Wiener processes. The $\dfrac{\sigma^2}{R}$ terms arise from It\^{o}'s lemma \cite{daffertshofer1998effects,gardiner1985handbook}. By using different values of parameters ($\beta,\delta$), we can control how fast the trajectory leaves the stable periodic orbit and converges to the stable equilibrium (having assumed that we start from a stable periodic solution). It is clear, from both the model formulation as well as \cref{modelsimulation} depicting a model simulation illustrating the transition from periodic behavior (stepping) to equilibrium (freezing) that this transition occurs with a uniform probability distribution in $\theta$ (i.e independent of the phase). \Cref{delembedding}{b} suggests that this may not be necessarily true in the case of freezing of gait in Parkinson's disease. 
\section{Robust location of escape from oscillations}
\label{sec:method}
\subsection{Hilbert Transform embedding of  stepping data time series} 

\begin{figure}[ht]
	\centering
		\setlength{\abovecaptionskip}{0.cm}
		\setlength{\belowcaptionskip}{-0.cm}
		\includegraphics[width=\textwidth]{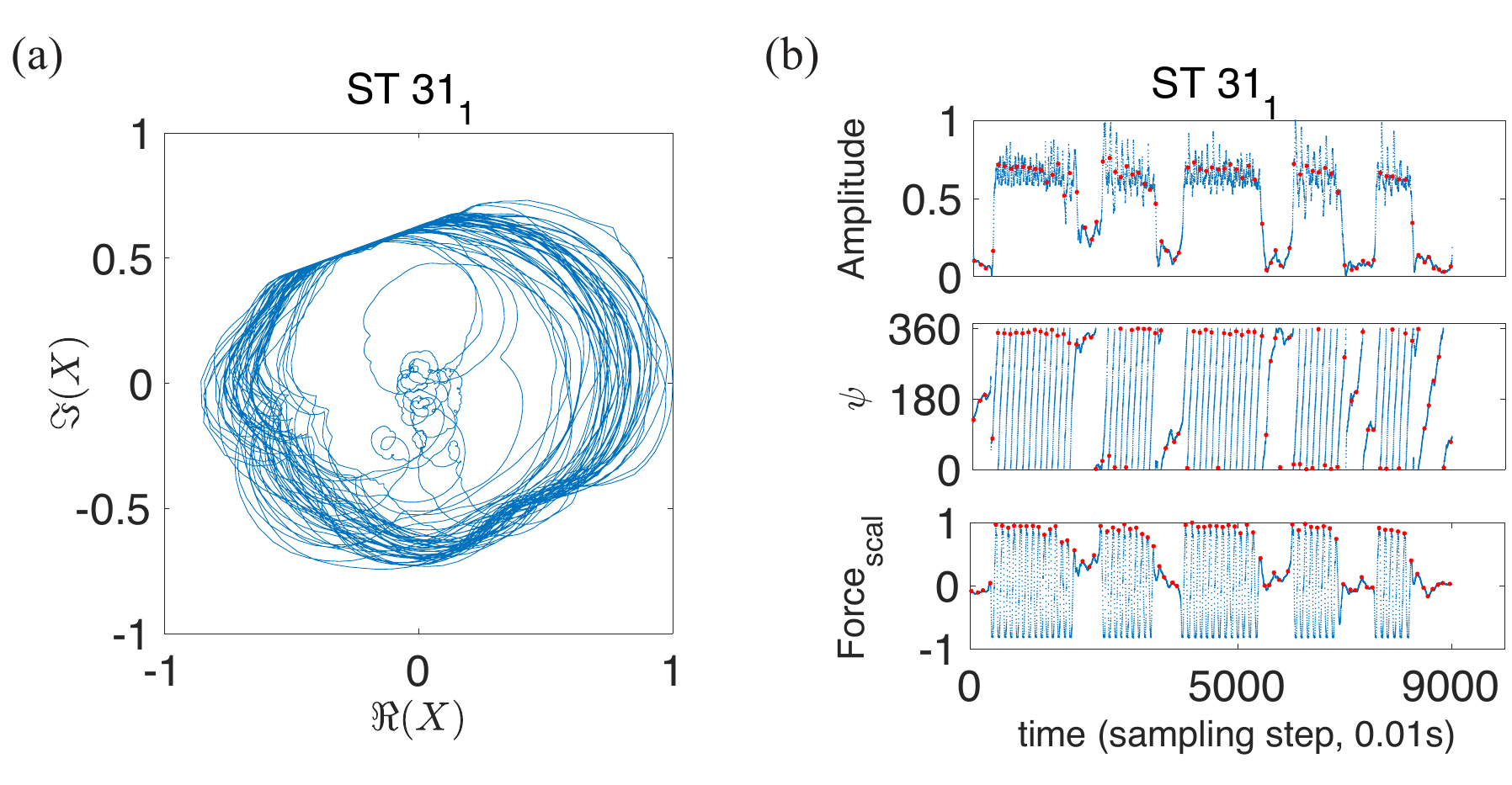}
		\caption{Panel(a): Embedding with scaled Hilbert Transform applied to stepping data from \cref{p1} (see also right, bottom panel). Panel (b): Amplitude and Phase ($\psi$)  of embedding. The red dots are locally maximum forces and their corresponding amplitude and phase.}
		\label{ht}
	\end{figure}

The nearly planar oscillatory
behavior during regular stepping in the stepping-in-place experimental data, visible in \Cref{delembedding}(a,b), makes the Hilbert Transform embedding
a natural choice when expressing stepping as an oscillator. More
precisely \cite{FELDMAN2011735}, for
a real scalar signal $x(t)$  on an interval of length
$T=2\pi/\omega$, expressed in terms of Fourier coefficients $x_k$, the
Hilbert Transform $\tilde{x}(t)$ has the Fourier coefficients
$-\mathrm{i}x_k$:
\begin{align*}
  x(t)&=x_0+\sum_{k=1}^\infty\left[x_k\mathrm{e}^{k\omega\mathrm{i}}+\bar x_k\mathrm{e}^{-k\omega\mathrm{i}}\right]\mbox{,}&
  \tilde{x}(t)&=\sum_{k=1}^\infty\left[-\mathrm{i}x_k\mathrm{e}^{k\omega\mathrm{i}}+\mathrm{i}\overline{x}_k\mathrm{e}^{-k\omega\mathrm{i}}\right]\mbox{,}
\end{align*}
such that the embedding is a complex scalar signal:
\begin{align*}
  X_\mathrm{hilbert}(t)=x(t)+\mathrm{i}\tilde{x}(t)\in\mathbb{C}\mbox{.}
\end{align*}
For discrete finite time series of stepping data the corresponding
discrete Fourier Transform (based on FFT) is used.  The signal
$X_\mathrm{hilbert}(t)$ oscillates around a non-zero mean with an
amplitude that is determined by the dimensions of the force
measurements. To obtain non-dimensionalized quantities for further
analysis we use the re-scaled signal $X(t)$:
\begin{align}\label{scaling}
  X(t)=\frac{X_\mathrm{hilbert}(t)-\operatorname{mean}(X_\mathrm{hilbert})}{\max|X_\mathrm{hilbert}(t)-\operatorname{mean}(X_\mathrm{hilbert})|}\mbox{.}
\end{align}
Using the representation of the signal $X(t)$ in its trigonometric or exponential form one can determine its instantaneous amplitude $R$ and phase $\psi$:
$$X(t) = \left| X(t)\right|\left[ \cos\psi(t)+\mathrm{i}\sin\psi(t) \right] = R(t)\mathrm{e}^{\mathrm{i}\psi(t)}\mbox{.}$$
In \cref{ht}(a) we show the Hilbert Transform embedding of the stepping time profiles  in the complex plane. 
The part of the approximate phase portrait following an approximate ellipse corresponds to regular stepping  episodes and the excursions toward the origin correspond to freezing episodes. \Cref{ht}(b) shows the instantaneous amplitude $R(t)$ and phase $\psi(t)$ of $X(t)$, and the original normalized data, respectively. The red dots in the \Cref{ht}(b) mark the force maxima (as determined by graph in the bottom panel) of each oscillation also in time-amplitude and time-phase plots. We can see that these peaks occur near $\psi=2\pi$, corresponds to the phase at which the left foot of the subject (in this case ST31) reaches the ground.	

  	\subsection{Selection of transition intervals for individual freezing events}\label{sec3}
\begin{figure}[ht]
		\centering
		\setlength{\abovecaptionskip}{0.cm}
		\setlength{\belowcaptionskip}{-0.cm}
		\includegraphics[width=\textwidth]{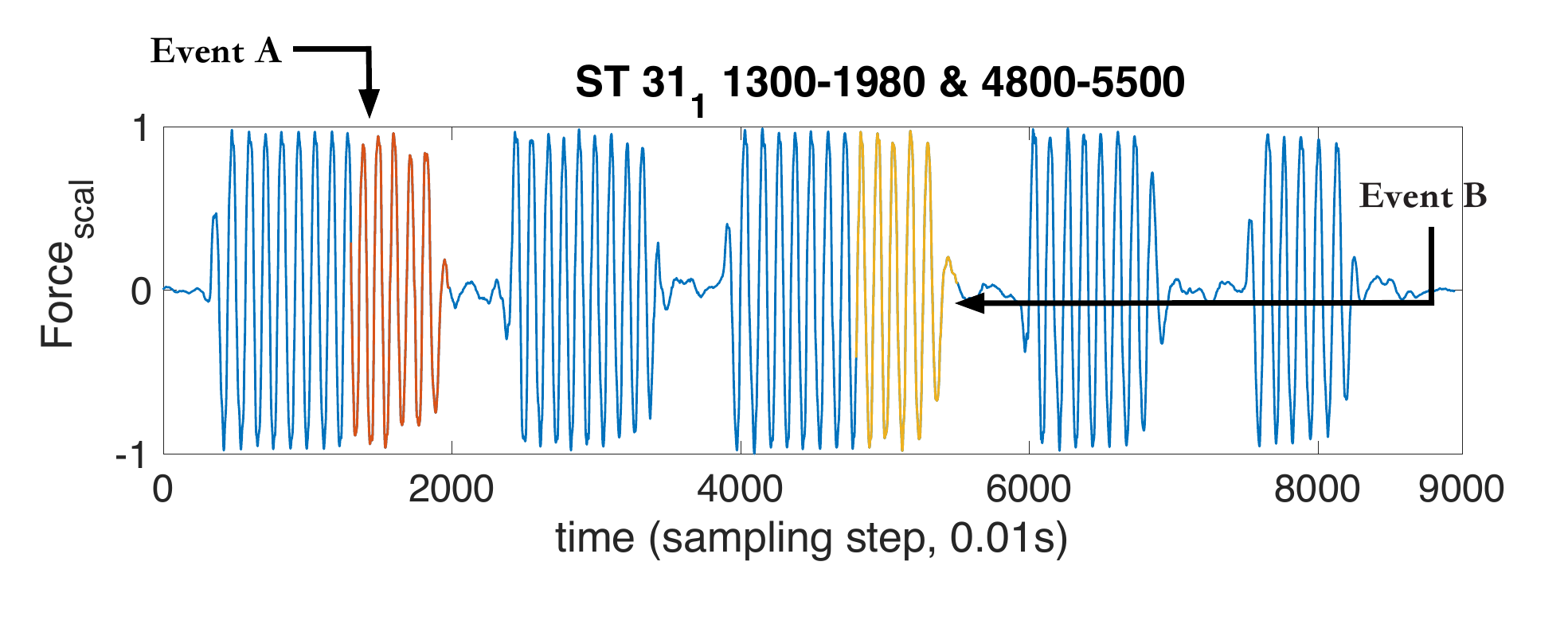}
		\caption{Sample events extracted from the original time series.  Event A represents a portion of time series from ST31, data set 1 in subscript and time interval [1300,1980]. Event B represents another portion of time series from ST31, data set 1 in subscript and time interval [4800,5500].}
		\label{new_events}
	\end{figure}

As we aim to study transitions from stepping to freezing we identify, for each freezing event $k$ labelled by a domain expert, a \emph{transition interval} $[t_{\mathrm{start},k},t_{\mathrm{end},k}]$ that contains this transition. \cref{new_events} highlights two examples (labelled $k=\{A,B\}$) from the time series shown in \cref{p1} (ST31, data set 1). The interval boundaries are  chosen such that $[t_{\mathrm{start},k},t_{\mathrm{end},k}]$ contains at least $4$ stepping oscillations and $t_{\mathrm{end},k}$ is inside the part of the time series identified as freezing.
For example, \cref{new_events} illustrates the highlighted events $A$ and $B$ extracted from ST31 (data set 1). Their transition intervals contain stepping and the beginning of the respective freezing episodes, in this case,
\begin{align*}
  [t_{\mathrm{start},A},t_{\mathrm{end},A}]&=[1300,1980]\mbox{,}&
  [t_{\mathrm{start},B},t_{\mathrm{end},B}]&=[4800,5500]\mbox{.}
\end{align*}
As the choice of $t_{\mathrm{start},k}$ and $t_{\mathrm{end},k}$ is arbitrary, we will for all results below determine how they depend on the choice of $t_{\mathrm{start},k}$ and $t_{\mathrm{end},k}$ for each event $k$.

\begin{figure}[ht]
		\centering
		\setlength{\abovecaptionskip}{0.cm}
		\setlength{\belowcaptionskip}{-0.cm}
		\includegraphics[width=\textwidth]{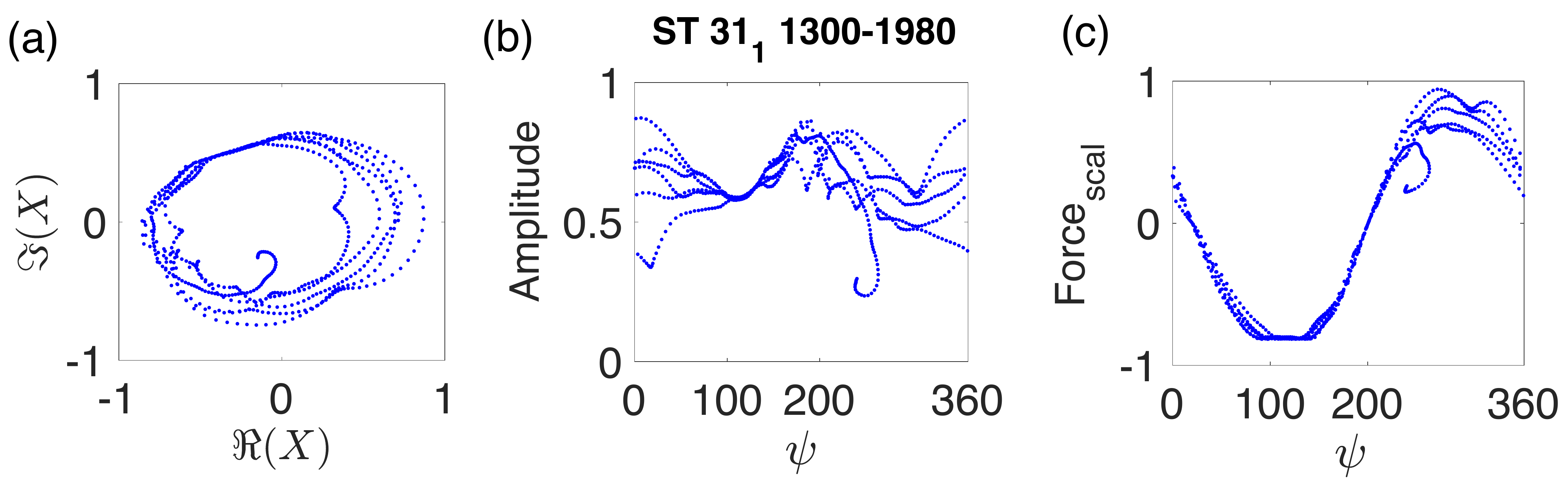}
		\caption{Embedded time series shown in different coordinate systems. Panel (a): embedded trajectory with Hilbert Transform. Panel (b): $x$-axis is angle $\psi$ in degrees, $y$-axis is amplitude. Panel (c): $x$-axis is angle $\psi$ in degrees, $y$-axis is scaled force. }
		\label{hilbert31info}
	\end{figure}

\cref{hilbert31info}(a) shows the data with Hilbert Transform for a transition interval that corresponds to freezing events extracted manually from the data (ST31, data set 1, time interval [1300,1980]). According to our convention, each panel's label shows the subject number with subscript indicating the number of the data set (we note that each subject has repeated the experiment several times and hence we have more than one data recordings set for each participant), followed by the pair of transition interval boundaries $[t_{\mathrm{start},k},t_{\mathrm{end},k}]$. As expected for transition intervals corresponding to a single freezing episode, the embedded trajectories initially follow an ellipse during regular stepping before they approach the area near the origin of the complex plane (only once, in contrast to \cref{ht}(a), since each transition interval contains exactly one freezing episode and respectively one transition). 

 \cref{hilbert31info}(b) and \cref{hilbert31info}(c) show the same set of embedded trajectories in two further coordinate systems, polar coordinates ($x$-axis is phase $\psi$ in degrees, $y$-axis is amplitude) in  \cref{hilbert31info}(b), and Cartesian coordinates $x$-axis is angle $\psi$ in degrees, $y$-axis is scaled force: $\mathrm{Force}_\mathrm{scal}$ given in \eqref{force_scaling}) in  \cref{hilbert31info}(c).

The representative transition in  \cref{hilbert31info} appears at first sight well described by the spontaneous transitions to the origin as they occur in the null model established by the Hopf normal form \eqref{eq:HNF}, up to a coordinate transformation that maps the ellipse followed during regular stepping onto the unit circle. However, we hypothesize a breaking of rotational phase symmetry, which would be a qualitative and practically relevant  difference to the Hopf normal form. The following sections will describe a method to identify the location of the transition in the embedded phase plane uniquely. Then we will apply this method to all freezing transitions in the available data sets to investigate our hypothesis provisionally. 

    		\subsection {Motivation for discrete-time discrete-space Markov chain}\label{sec_rw}
\Cref{hilbert_polar1,amplitude_degree1}  in \cref{sec:more:events} show a survey of Hilbert embeddings for freezing transition intervals from patient data. \Cref{hilbert_polar1} uses cartesian coordinates $X$ and \cref{amplitude_degree1} uses polar coordinates $(\psi,R)$. Especially, \cref{amplitude_degree1} highlights that fixing a threshold (for, e.g., amplitude $R$ in polar coordinates as $R_\mathrm{th}$) and defining the time of transition (for example) as the first crossing of this threshold (e.g., first time $t_\mathrm{th}\in[t_{\mathrm{start},k},t_{\mathrm{end},k}]$ when $R(t_\mathrm{th})\leq R_\mathrm{th}$) will introduce an extreme dependence of the timing $t_\mathrm{th}$ and the angle coordinate $\psi(t_\mathrm{th})$ on this threshold value. Furthermore, the threshold value will have to be adjusted for different events (possibly even between events for the same subject and data set), making the collection of generalizable statistics across subjects and events impossible and reinforcing the need for a subject-specific or personalized approach. In contrast to the sensitivity of timing and angle of the transition to the threshold,  the \emph{detection} of the presence of a freezing event with the help of thresholds is robust. As mentioned in \cref{sec:data:intro}, freezing is uniquely determined by the absence of zero-force periods for longer than the identified stepping frequency. 

We exploit the unique identifiability of the presence of freezing by constructing a discrete-time discrete-space Markov chain and its associated transition matrix for the variable $X$ obtained from the Hilbert Transform embedding. This Markov chain divides the phase space into a transition set and an absorbing set (see below for precise definitions). The starting time into the freezing event can then be uniquely defined as the transition between these two subsets of the phase space. The red crosses in \cref{hilbert_polar1,amplitude_degree1} show the transition points resulting from the procedure described below. They are clearly not determined by an amplitude coordinate, but they are independent of discretization parameters for our Markov chain  (see \cref{metanalysis31}) and independent of the selected transition interval $[t_{\mathrm{start},k},t_{\mathrm{end},k}]$ (see \cref{convar_analysis}). 

The Markov chain also allows us to generate ``surrogates'': time series $X_\mathrm{mc}(t)$ in the complex plane that share properties with the embedded time series $X(t)$ of the  data shown in \cref{hilbert31info}(a) for each freezing event and each subject.

\subsection{Subdivision of the complex plane along polar coordinates}\label{sec:subdivison}
\begin{figure}[ht]
          \centering
          \setlength{\abovecaptionskip}{0.cm}
          \setlength{\belowcaptionskip}{-0.cm}
          \includegraphics[scale=0.4]{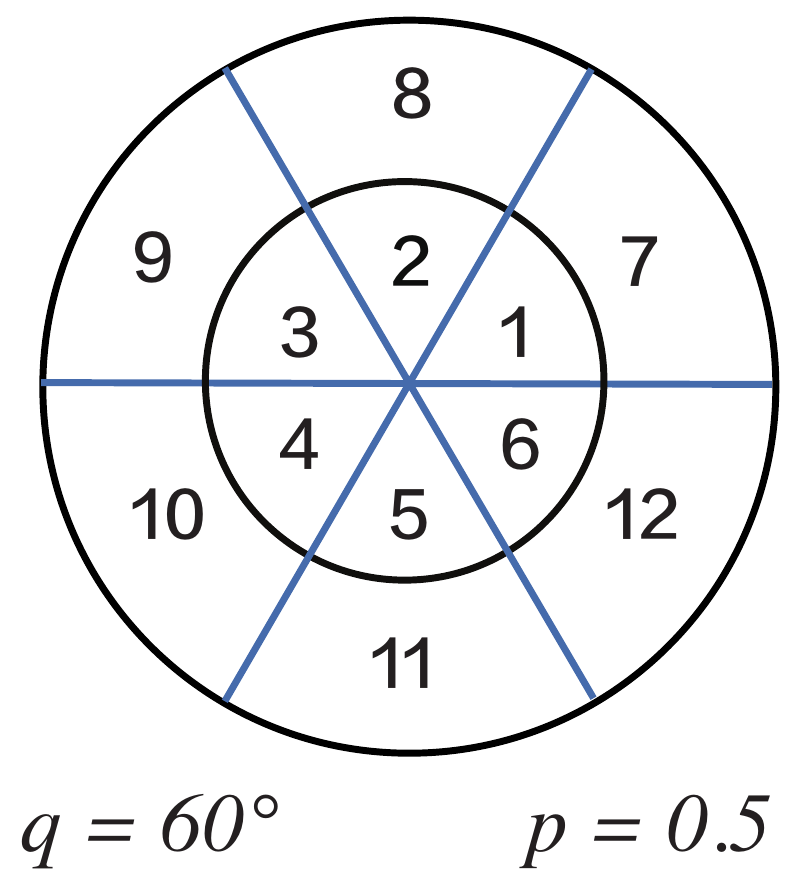}
          \caption{Subdivision and discretization of the unit circle in the complex plane. The numbering is anticlockwise and from the pole to the edge of the space with $P=2$, $Q=6$ and, hence, $p=0.5$ and $q=60^\circ$. Note that typical tested discretizations are much finer: \cref{metanalysis31} tests the range $q\in\{3^\circ,5^\circ,10^\circ,15^\circ,20^\circ,30^\circ\},p\in\{0.05,0.1,0.15,0.2\}$.}
          \label{subdivide}
        \end{figure}

We subdivide the complex plane (the embedding space after Hilbert Transform) into boxes that are rectangular in polar coordinates as illustrated in \cref{subdivide}. 
Thus, the box boundaries are aligned with the radial direction and angular direction, respectively.
\begin{itemize}
\item By the scaling of our data the amplitude of $X(t)$ is always in $(0,1)$, such that we subdivide in the radial direction into $P$ annuli of equal radial thickness $p=1/P$. 
  We enumerate starting from the origin. 
\item We subdivide in the angular direction anticlockwise, starting from positive $x$-axis into $Q$ cones of equal size $q=360^\circ/Q$.
\end{itemize} 
This results in $P\times Q$ discrete states in total, corresponding to the near rectangular boxes described above.
A box $B_{k,\ell}\subset\mathbb{C}$ for $(k,\ell)\in\{1,\ldots,P\}\times\{1,\ldots,Q\}$ is then 
\begin{align*}
  B_{k,\ell}=\left\{X\in\mathbb{C}:R(X)\in\left(\frac{k-1}{P},\frac{k}{P}\right],
  \psi(X)\in\left[\frac{2\pi(\ell-1)}{Q},\frac{2\pi\ell}{Q}\right)\right\}\mbox{,}
\end{align*}
where we use $R(X)\in(0,1]$ for the amplitude and
$\psi(X)\in[0,360^\circ)$ for the argument of a point
$X\in\mathbb{C}$. We enumerate the boxes in angle-first order such
that the box with radius-angle index $(k,\ell)$ is at position
$i=\ell+Q(k-1)$. The above subdivision defines index maps
$\ir:\mathbb{C}\to\{1,\ldots,P\}$, 
$\ia:\mathbb{C}\to\{1,\ldots,Q\}$ and $\ind:\mathbb{C}\to\{1,\ldots,PQ\}$ as follows:
\begin{align}
  \ir(X)&=\ceil(R(X)P) 
  \label{P_i}\\
  \ia(X)&=\floor(\psi(X)Q/(2\pi))+1  
  \label{Q_i}\\
  \ind(X)&=(\ir(X)-1)Q+\ia(X)\mbox{.}
  \label{X_i}
\end{align}
For box number $i$ we may recover the corresponding annulus $k$ and cone $\ell$ via
\begin{align}
  \ell&=i-Q(\floor(i/Q)-1)\mbox{,}&
  k&=\floor(i/Q)\mbox{,} \label{radang:frombox}
\end{align}
such that
\begin{align*}
  \ia(X)&=\ind(X)-Q\times (\floor(\ind(X)/Q)-1)\mbox{,}&
  \ir(X)&=\floor(\ind(X)/Q)\mbox{}  
\end{align*}
for every $X\in\mathbb{C}$.

        In principle, we may consider the set $\{1,\ldots,PQ\}$ as the state space of our discrete-time discrete-space Markov chain. However, each box index corresponds to an approximate location inside the complex unit circle, given by its mid point with respect to radius and angle:
        \begin{align}\label{def:Xc}
          X_{\mathrm{c},i}=\frac{k-1/2}{P}\exp\left(\frac{2\pi\mathrm{i}(\ell-1/2)}{Q}\right)\mbox{,}
        \end{align}
        where $k$ and $\ell$ are related to $i$ via \eqref{radang:frombox}. Thus, $X_\mathrm{c}$ is a vector in $\mathbb{C}^{PQ}$, and for each $X\in\mathbb{C}$ we can find the midpoint of the box in which $X$ is located by
        \begin{align*}
          X_\mathrm{c}(X)&=X_{\mathrm{c},\ind(X)}\mbox{,\quad with coordinates}\\
          R_\mathrm{c}(X)&=\frac{\ir(X)-1/2}{P}\\
          \psi_\mathrm{c}(X)&=\frac{2\pi\mathrm{i}(\ia(X)-1/2)}{Q}\mbox{.}
        \end{align*}
        We may use $X_\mathrm{c}$ as the state space of the
        discrete-space Markov chain instead of the index in
        $\{1,\ldots,PQ\}$.
	\subsection{An empirical Markov chain transition matrix for a single transition interval}
	For the discrete state space $X_\mathrm{c}\sim\{1,\ldots,PQ\}$, defined in \eqref{def:Xc}, we use the embedded trajectory $X(t)\in\mathbb{C}$ with $t$ in a single transition interval $[t_\mathrm{start},t_\mathrm{end}]$ to construct a provisional transition matrix $A_\mathrm{c}$ for probability distributions $P:\{1,\ldots,PQ\}\to[0,1]$: 
    \begin{align}\label{gen:transmatrix}
	P_{n+1} = P_n A_\mathrm{c}
      \end{align}
	For a given embedded trajectory $X(t)$ with $t\in[t_\mathrm{start},t_\mathrm{end}]$ (recall that $t$ were integers) we define the transition count $\mathrm{CT}\in\mathbb{Z}^{PQ\times PQ}$ as follows:
        \begin{align}
          \label{def:transitioncount}
          \mathrm{CT}_{ij}=\left|\left\{t\in[t_\mathrm{start},t_\mathrm{end}-1]:\ind(X(t))=i\mbox{\ and\ }\ind(X(t+1))=j\right\}\right|
        \end{align}
        for $i,j\leq PQ$ (we use the notation $|S|$ for the number of elements in a set $S$ in \eqref{def:transitioncount}). The provisional transition matrix $A$ generated by $X(t)$ is then:
        \begin{align}
          \label{def:transitionmatrix}
          \{A_\mathrm{c}\}_{ij}=
          \begin{cases}
            \cfrac{\mathrm{CT}_{ij}}{\sum_{k=1}^{PQ}\mathrm{CT}_{ik}}&\mbox{if $\sum_{k=1}^{PQ}\mathrm{CT}_{ik}>0$,}\\[2ex]
            0&\mbox{otherwise}
          \end{cases}
        \end{align}

        The entry $ \{A_\mathrm{c}\}_{ij}$ is then the empirical probability (as determined by counting) that $X(t+1)$ is in state $j$ if $X(t)$ is in state $i$. By construction the row sums of $A$ satisfy $\sum_{j} \{A_\mathrm{c}\}_{ij} = 1$ or $\sum_{j} \{A_\mathrm{c}\}_{ij} = 0$. We call
        \begin{align}
          X_\mathrm{emp}&=\{X_{\mathrm{c},i}:\sum_{k=1}^{PQ}\mathrm{CT}_{ik}>0\}\mbox{\quad (in complex coordinates),}\\
          \mathcal{J}_\mathrm{emp}&=\{i:\sum_{k=1}^{PQ}\mathrm{CT}_{ik}>0\}\subseteq\{1,\ldots,PQ\}\mbox{\quad (in integer coordinates),}&
        \end{align}
        \begin{align*}
            n_\mathrm{emp} = |\mathcal{J}_\mathrm{emp}|
        \end{align*}
        the state space with empirical support. 
	
	\begin{figure}[ht]
\centering
					\setlength{\abovecaptionskip}{0.cm}
			\setlength{\belowcaptionskip}{-0.cm}
			\includegraphics[width=\textwidth]{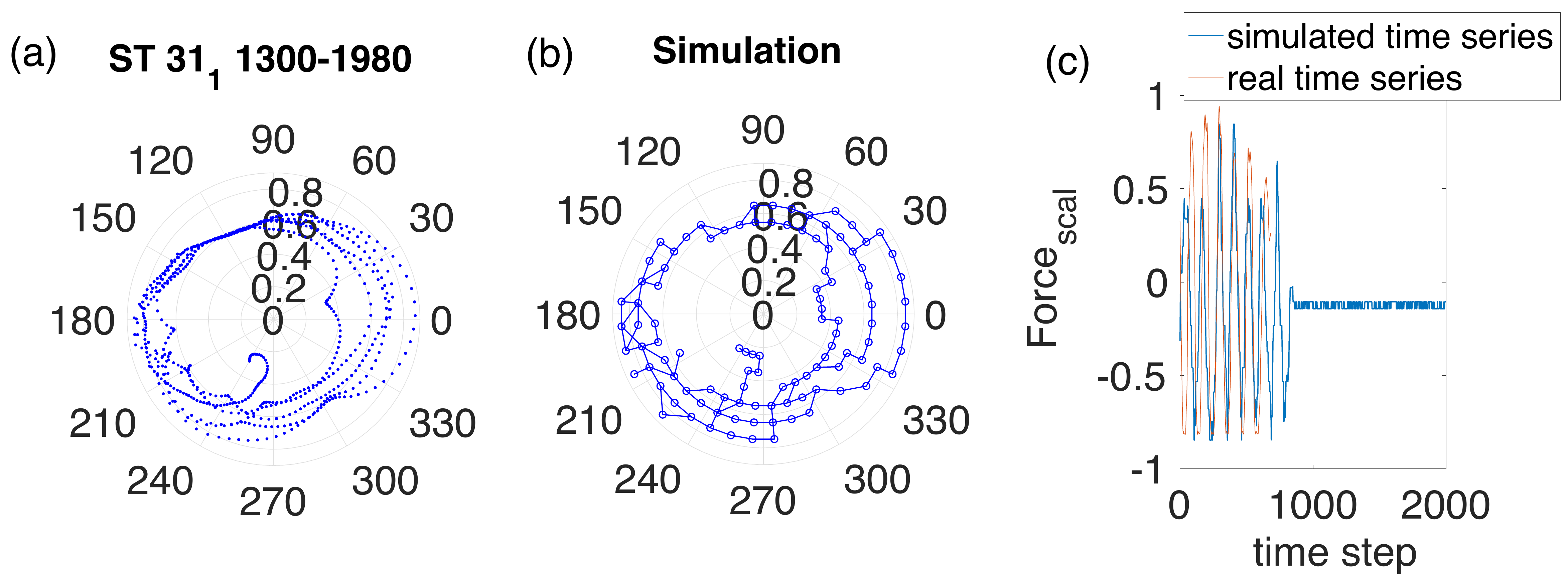}
	\caption{Embedded trajectory and surrogate time series for ST31, data set $1$, transition interval $[1300,1980]$ and discretization parameters $(p,q)=(0.1,5^\circ)$. Panel (a) shows the embedded trajectory $X(t)$. Panel (b) demonstrates a single surrogate time series while Panel (c) shows the projection onto the real part, the coordinate of the embedding corresponding to the scaled forced   data (blue: surrogate, red: data). }
		\label{random walk}
	\end{figure}
	
        We construct an empirical transition matrix $A$ which is restricted to $\mathcal{J}_\mathrm{emp}$. The transition matrix $A$ with state space $X_\mathrm{emp}$ allows us to either generate surrogate time series by starting from a state $X_{\mathrm{emp},i}$ at time $0$ and then generating the state at time $n+1$ by performing a multi-outcome Bernoulli trial with outcomes in $X_\mathrm{emp}$ and probabilities $A_{i,\mathcal{J}_\mathrm{emp}}$, and then changing the state at time $n+1$ to the outcome of the trial. 
        \Cref{random walk}(b) depicts such a surrogate time series. \Cref{random walk}(c) shows the same time series, superimposing it onto the real component of the embedded trajectory, which is the coordinate of the original signal.

\section{Properties of the constructed Markov chain}
\label{sec:mc:prop}
	\subsection{Partition of state space into  classes}
        A property of our empirical data $X(t)$ and the selection of transition intervals is that they generate Markov chain transition matrices $A$ with a distinct subdivision into communicating classes. We recall the following two standard definitions, where we let $\mathbf{u}_i$ be the $i$th unit row vector ($u_{ij}=1$ if $i=j$, $0$ otherwise). 
        \begin{defn}[accessible]
		A state $j$ is said to be accessible from a state $i$ (written $i\rightarrow j$) if there exists $n>0$ such that $\mathbf{u}_iA^n\mathbf{u}_j^\mathsf{T}>0$ \cite{markov1971extension}.
	\end{defn}
	
	\begin{defn}[Communicating class]
		A state $i$ is said to communicate with state $j$ (written $i\leftrightarrow j$) if both $i \rightarrow j$ and $j \rightarrow i$. The relation $\leftrightarrow$ is an equivalence relation. A communicating class is an equivalence class of $\leftrightarrow$, so, it is a maximal set of states $G$ such that every pair of states in $G$ communicates with each other \cite{markov1971extension}.
	\end{defn}
	Communicating classes only depend on the positivity of entries $A_{ij}$ such that we may introduce the indicator matrix  $Z$ and the reachability matrix $Z_\infty$:
        \begin{align*}
          Z_{ij} &=
          \begin{cases}
            1& \mbox{if $A_{ij}>0$,}\\
            0& \mbox{otherwise,}
          \end{cases}\mbox{,}&
          Z_{\infty,ij} &=
          \begin{cases}
            1& \mbox{if $A_{ij}^n>0$ for some $n>0$,}\\
            0& \mbox{otherwise.}
          \end{cases}
        \end{align*}
        If we  consider states $i$ and $j$ connected when $A_{ij}>0$, then the transition matrix $A$ induces a graph, for which $Z$ is the adjacency matrix. Thus, we may use graph theoretic methods to determine the reachability matrix $Z_\infty$ (in Matlab the routine \texttt{transclosure}).
From $Z_\infty$ we may construct the
bidirectional reachability matrix $Z_\mathrm{bi}=\min(Z_\infty,Z_\infty^\mathsf{T})$, that is,
\begin{align*}
  Z_{\mathrm{bi},ij} =
  \begin{cases}
    1 &\mbox{if $Z_{\infty,ij}=1$ and $Z_{\infty,ji}=1$,}\\
       0&\mbox{otherwise.}
      \end{cases}
\end{align*}
Two states $i$ and $j$ are in the same communicating class if the rows $Z_{\mathrm{bi},i,(\cdot)}$ and $Z_{\mathrm{bi},j,(\cdot)}$ are
identical. Let the \emph{class indicator matrix}
$Z_\mathrm{cc}\in\{0,1\}^{n_\mathrm{cc}\times n_\mathrm{emp}}$ be the matrix consisting
of the $n_\mathrm{cc}$ unique rows of $Z_\mathrm{bi}$ (typically $n_\mathrm{cc}\ll n_\mathrm{emp}$), then each row of $Z_\mathrm{cc}$
corresponds to a communicating class of transition matrix $A$,
and state $i$ is in communicating class $k$ for $A$ if $Z_{\mathrm{cc},k,i}=1$. Communicating classes are partially ordered: we write that class
\begin{align*}
    k_1\rightarrow k_2
\end{align*}
 (with $k_1\neq k_2$) if states from class $k_2$
are reachable from states in class $k_1$ (if reachability is true
for one pair of states in classes $k_1$ and $k_2$ it is true for all pairs of states).

\subsection{Communicating classes of empirical transition
matrices --- stepping class, transition set and absorbing set} The transition intervals $[t_\mathrm{start},t_\mathrm{end}]$ in Section \ref{sec3} are
chosen such that the time series $x(t)$ starts in the stepping
regime (large scaled force oscillations), stays there for most
of the time and ends in a freezing episode. Thus, we expect
the states of the Markov chain to fall into several
communicating classes (see \cref{random walk}(a)). \Cref{discritise size} studies systematically
the dependence of the number of communicating classes, their ordering and their geometric properties (e.g., phase angles $\psi_\mathrm{c}$ of the boxes where transition between classes is most likely, see Section~\ref{sec:phase:prediction}) on discretization
parameters $P$ and $Q$ and transition interval
$[t_\mathrm{start},t_\mathrm{end}]$. Our results in
\cref{discritise size} show that this dependence is weak. Thus communicating classes are a suitable object for studying the qualitative properties of our time series data. We expect  communicating classes of our data from
class indicator matrix $Z_\mathrm{cc}$ of the following types and with the following properties.
\begin{itemize}
\item \textbf{(Stepping)} We expect one large communicating class corresponding to a row in the class indicator matrix $Z_\mathrm{cc}$ with almost all the values with empirical support equal to $1$. We denote this communicating class as the \emph{stepping
class}, naming it $F_\mathrm{step}$ and denoting its row index in $Z_\mathrm{cc}$ by $k_\mathrm{step}$. 
\item\textbf{(Ergodicity)} By definition of a communicating
class, all states $j$ in the stepping class $F_\mathrm{step}$ are positively recurrent. We also assume that they are aperiodic. This is an assumption that the
discretization box sizes in the Hilbert Transform plane
$\mathbb{C}$ should not be too small compared to the
sampling time step: at least some boxes from the
discretization should contain several sampling time steps
from the sampled time series such that the empirical
probability of staying in the box during a time step is
non-zero. With this aperiodicity assumption the stepping
class is ergodic when one considers the conditional
transition probabilities $P_{\mathrm{step},i,j}$  under the condition that the
Markov chain stays in the stepping class:
\begin{align}\label{fstep:ergodic}
 P_{\mathrm{step},i,j}=P(X_{n+1}=i\,|\,X_n=j\mbox{\ and\ }X_{n+1}\in
 F_\mathrm{step})\mbox{\ for $i,j\in F_\mathrm{step}$.}
\end{align}
\item \textbf{(Ordering relative to stepping class)} We also expect that all other classes can be related to $F_\mathrm{step}$ through the partial ordering (because $F_\mathrm{step}$ is large). That is, each  class $k$ satisfies $k\rightarrow k_\mathrm{step}$ or $k_\mathrm{step}\rightarrow k$.		
\item \textbf{(Transition set)} We collect the states in all communicating classes $k$ with $k\rightarrow k_\mathrm{step}$ (so coming before $F_\mathrm{step}$ in the partial ordering) in the so-called initial set $F_0$. The \emph{transition set} is the union of stepping class $F_\mathrm{step}$ and initial set $F_0$, $F =F_0\cup F_\mathrm{step}$.
\item \textbf{(Absorbing set)} We collect all communicating classes $k$ with $k_\mathrm{step}\rightarrow k$ (so coming after $F_\mathrm{step}$ in the partial ordering) in the so-called absorbing set $E$. We expect the absorbing set $E$ to be non-empty.
\end{itemize}

		\begin{figure}[ht]
		\centering
		\setlength{\abovecaptionskip}{0.cm}
		\setlength{\belowcaptionskip}{-0.cm}
		\includegraphics[width=\textwidth]{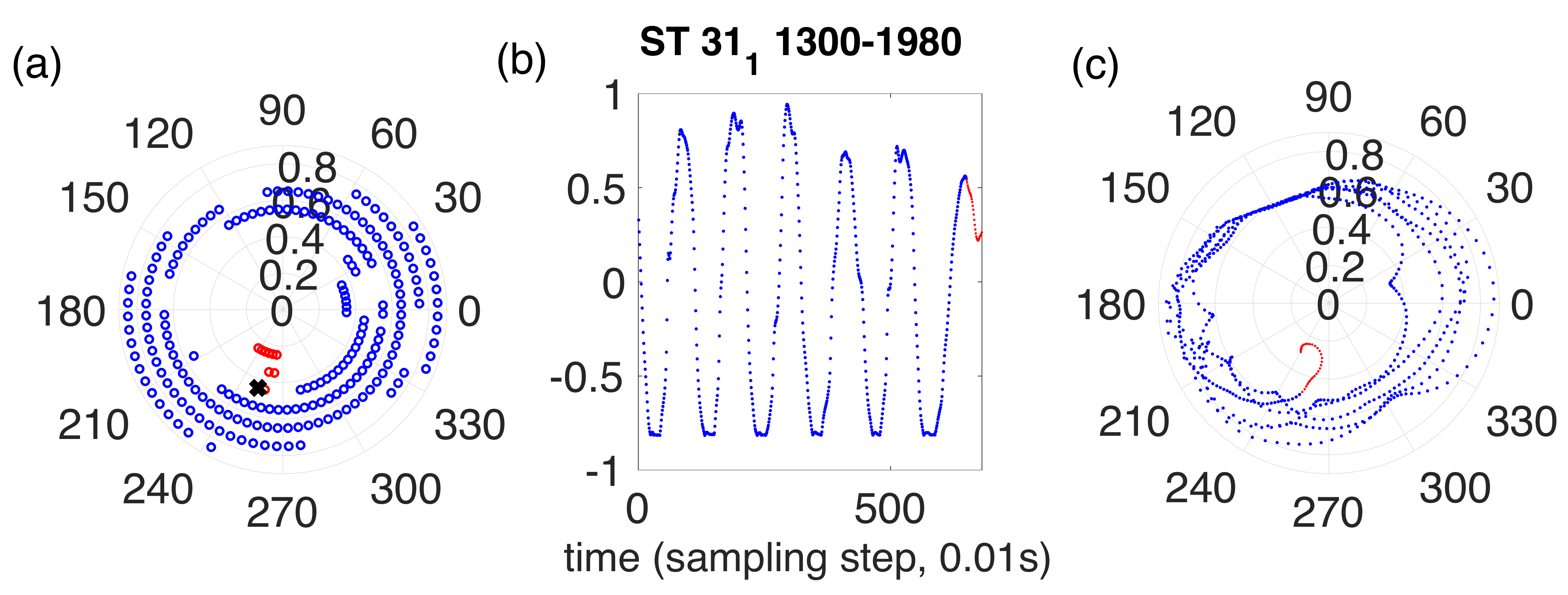}
		\caption{ Panel (a): empirical support for Markov chain for time series ST31 and decomposition into transition set $F$ (blue circles) and absorbing set $E$ (red circles) with discretization parameters $(p,q)=(0.1,5^\circ)$. The black cross is the ``first'' state, $X_\mathrm{tr}$, in the absorbing set $E$ (see \Cref{sec:transition}). Panel (b):  underlying time series for panel (a), color coding the sampling points according to their location in transition set (blue) or absorbing set (red). Panel (c): same time series and color coding in phase plane obtained by Hilbert Transform.}
		\label{tran31}
        \end{figure}   

\cref{tran31} shows an example for a transition matrix generated from ST31, data set 1, transition interval $[1300,1980]$, and discretization parameters $p=0.5$,
$q=5^\circ$. Panel (a) shows all boxes with empirical support in        the form of circles at their centers $X_{\mathrm{emp},i}$. The red circles are in boxes that belong to the absorbing set $E$,
and the blue circles are in boxes that belong to transition set $F$. We observe that there is a unique state $i$ in the absorbing set with $A_{ij}>0$ for a state $j$ in the transition set and highlight this point with a black cross. \Cref{tran31}(b) and (c) show the original time series $x(t)$ (unscaled force measurements) and its embedding in the unit circle of the complex plane. We color the points of the time series according to the class membership of the box they are in: red marking indicates that $\ind(X(t))$  is in the absorbing set $E$. Blue marking indicates that $\ind(X(t))$ is in the transition set $F$.
\subsection{Mean escape time from the transition
          set} \label{escape} Suppose we can decompose the Markov
        chain given by transition matrix $A$ into a transition set $F$
        with $m$ states and an absorbing set $E$ with $k=n-m$
        states. Let us denote the transition matrix, restricted to the
        transition set $F$ by $A_F\in\mathbb{R}^{m\times m}$. Let
        $\mathbf{s}\in\mathbb{R}^{1\times m}$ be an initial
        probability distribution for states in the transition set
        $F$. The probability that, starting from this distribution,
        one stays inside $F$ until at least time $t$ equals
        $\mathbf{s}A_F^t\mathbf{e}^\mathsf{T}$, where
        $\mathbf{e}=(1,\ldots,1)\in\mathbb{R}^{1\times
          m}$. Consequently, the probability of leaving the transition
        set exactly at time $t$ equals $\mathbf{s}A_F^{t-1}\mathbf{e}^\mathsf{T}-\mathbf{s}A_F^t\mathbf{e}^\mathsf{T}$ (the event of survival to time $t$ is a subset of the event of survival to time $t-1$). Hence, the mean escape time from the transition set, denoted $\MET(\mathbf{s}$), starting from initial distribution $\mathbf{s}$ equals
	\begin{align}
	\text{MET($\mathbf{s}$)} & =
	\sum_{t=1}^{\infty} t \left[ \mathbf{s}A_F^{t-1}\mathbf{e}^\mathsf{T}-\mathbf{s}A_F^t\mathbf{e}^\mathsf{T}\right]= \mathbf{s}\sum_{t=1}^{\infty} t A_F^{t-1} \left(I-A_F\right)\mathbf{e}^\mathsf{T}\label{rankproblem}
	\\&=\mathbf{s} \left(I-A_F\right)^{-1}\mathbf{e}^\mathsf{T}\mbox{.}\label{metbymatrix} %
\end{align}
As the absorbing set $E$ comes after the transition set $F$ in
the partial ordering of the communicating classes, every left
eigenvector $\bm{\pi}$ and eigenvalue $\lambda$ of $A_F$
        satisfies
        $[\bm{\pi}A_F^k]_j=[\lambda^k\bm{\pi}]_j<\bm{\pi}_j$ for
        at least one index $j\leq m$ and a sufficiently large $k$, since states from outside the support
        of $\mathbf{s}$ will always be reachable in sufficiently many ($k$) steps. Hence,
        $|\lambda|<1$ for all eigenvalues of $A_F$, such that $I-A_F$ is
        invertible.
Thus, the identity
$\sum_{t=1}^{\infty}tA_F^{t-1} = (I-A_F)^{-2}$, applied for \eqref{metbymatrix}, is justified.
\subsection{Preferred transition states}
\label{sec:transition}
\begin{figure}[ht]
		\centering
		\setlength{\abovecaptionskip}{0.cm}
		\setlength{\belowcaptionskip}{-0.cm}
		\includegraphics[width=0.8\textwidth]{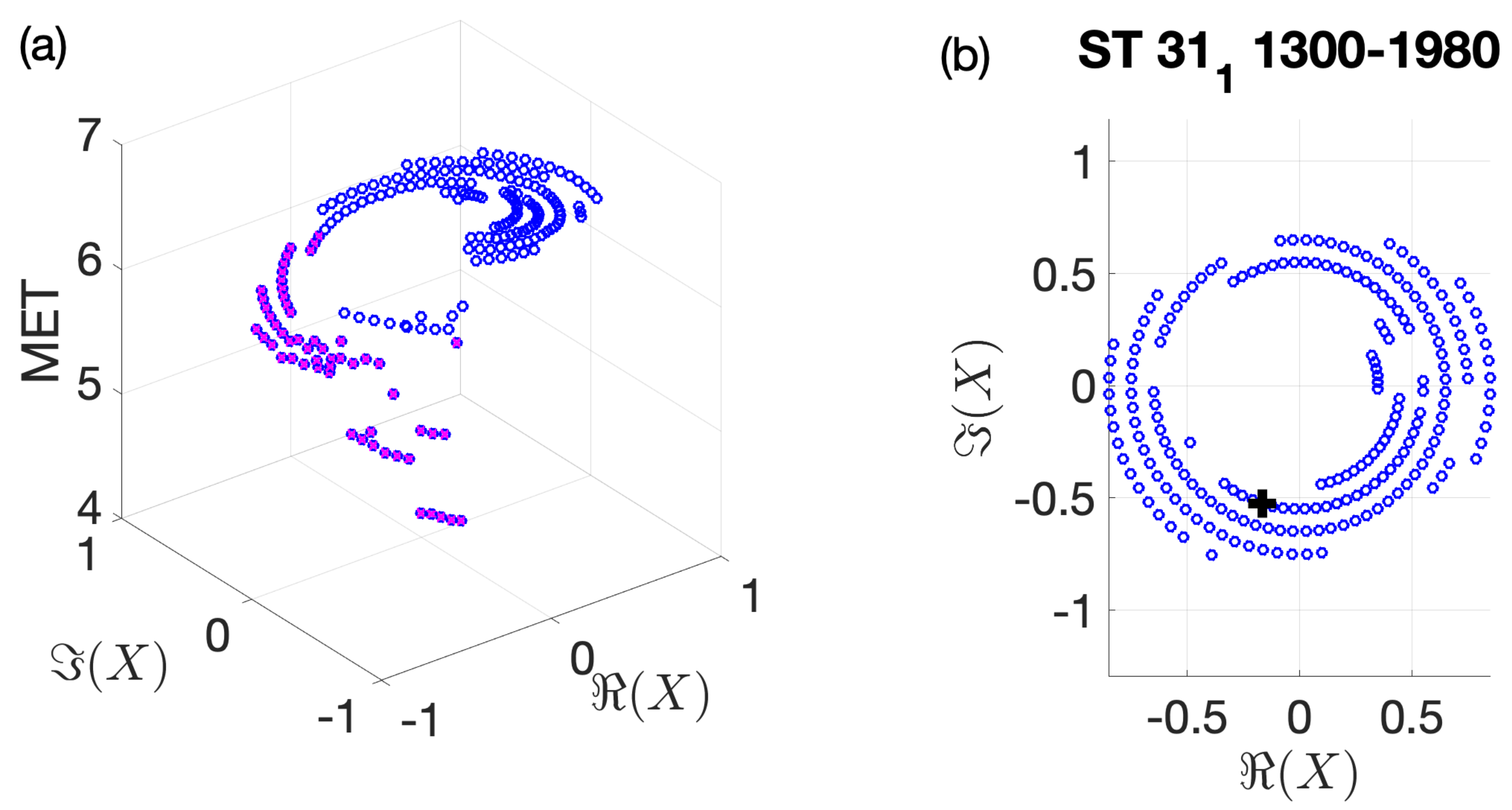}
		\caption{Panel (a): mean escape times from each state in transition set $F$ for ST31 data set 1, time interval $[1300,1980]$, discretization parameters $(p,q)=(0.1,5^\circ)$. States marked by magenta crosses have $\MET$ below $\MET_F$ given in \eqref{metbylambda1}. Panel (B): top view, where black '+' marks $X_{\min}$, the state with minimal mean escape time.}
		\label{metconclusion31}
	\end{figure}

\cref{metconclusion31} indicates the mean escape times for $i\in F$ and initial distributions equaling unit vectors,
\begin{align}\label{def:meti}
  \MET_i=\MET(\mathbf{u}_j)=\mathbf{u}_j(I-A_F)^{-1}\mathbf{e}^\mathsf{T}\mbox{,~where $i$ is the $j$th element of $F$,}
\end{align}
and $u_{j,j}=1$ and $u_{j,\ell}=0$ for $\ell\neq j$ ($\mathbf{u}_j\in\mathbb{R}^{1\times m}$). Thus, $\MET_i$ is the mean escape time when we are starting from a known state $i$ in the transition set $F$. \Cref{metconclusion31}(b) shows the top view and highlights the state $i_{\min}$ with minimal $\MET_i$ and its complex coordinate $X_\mathrm{min}$:
        \begin{align}
          \label{eq:imin}
          i_\mathrm{min}&=i\mbox{\ such that such that $\MET_i$ is minimal, and}\\
          \label{eq:xmin}
          X_\mathrm{min}&=X_{\mathrm{emp},i_\mathrm{min}}\mbox{,\quad}
          R_\mathrm{min}=R(X_{\mathrm{emp},i_\mathrm{min}})\mbox{,\quad}
          \psi_\mathrm{min}=\psi(X_{\mathrm{emp},i_\mathrm{min}})\mbox{.}
        \end{align}
The center of the box $X_\mathrm{min}$ is marked by a black '+' in \cref{metconclusion31}(b). The location of $X_\mathrm{min}$ is close to the (in this case unique) ``first'' state $X_\mathrm{tr}$ (with amplitude $R_\mathrm{tr}=R(X_\mathrm{tr})$, phase $\psi_\mathrm{tr}=\psi(X_\mathrm{tr})$ and index $i_\mathrm{tr}\in E$) in the absorbing set into which one may transition from $F$, which is highlighted by a black cross in \cref{tran31}(a). If $i_\mathrm{tr}$ is unique, then it is its own communicating class such that it comes indeed first in the partial ordering within $E$ and follows directly in the partial ordering after $F$. The boxes $X_\mathrm{tr}$ and $X_\mathrm{min}$ are naturally close together whenever a unique $i_\mathrm{tr}\in E$ exists in the sense of our partial ordering. While the mean transition times will clearly depend on the choice of boundaries for the transition interval $[t_\mathrm{start},t_\mathrm{end}]$, the transition set $F$, the absorbing set $E$ and $X_\mathrm{min}$ and $X_\mathrm{tr}$ may not (or only weakly) depend in $t_\mathrm{start}$ or $t_\mathrm{end}$. Similarly, discretization parameters $P$ and $Q$ may only weakly affect $X_\mathrm{min}$ and $X_\mathrm{tr}$. If these weak or non-dependences are true then a systematic collection of $X_\mathrm{min}$ (and $\psi_\mathrm{min}$) from subject data sets will be able to determine whether the phase invariance implicitly assumed in the generalized Hopf normal form model is a valid assumption, or if (subject dependent or general) phases exist during which subjects are particularly vulnerable to freezing. \cref{discritise size} studies the dependence of $X_\mathrm{min}$ on all method parameters in detail.

	\subsection{Relationship between $\MET_i$ and the mean escape time from transition set $F$}
	The largest eigenvalue $\lambda_1$ in modulus (if there exists a largest in modulus eigenvalue) for the transition matrix $A_F$ of the transition set $F$ is the escape rate for a particular distribution $\mathbf{s}_F$, corresponding to the left eigenvector of $A_F$ for $\lambda_1$. This distribution $\mathbf{s}_F$ is the so-called stationary distribution of the transition set $F$, that is, $\mathbf{s}_F$ is the limiting distribution of distributions $\mathbf{s}$,
	\begin{align*}
	    \mathbf{s}_N=\mathbf{s}A_F^N/(\mathbf{s}A_F^N\mathbf{e}^\mathsf{T})
	    \to\mathbf{s}_F,  
	\end{align*}
	after $N$ steps under the condition that the Markov chain stays in the transition set $F$ for at least $N$ steps, and $N\to\infty$. The eigenvalue $\lambda_1$ is called the mean survival probability (per time step) for $F$, such that the mean escape time for $F$ equals
	\begin{align}
	\MET_F:=\frac{1}{1-\lambda_1}=\MET\left(\mathbf{s}_N\right)+O((\lambda_\mathrm{dec}/\lambda_1)^N)\mbox{,} \label{metbylambda1}
	\end{align}
	where $\lambda_\mathrm{dec}$ is the modulus of the second-largest eigenvalue of $A_F$. The mean survival time $\MET_F$ in the Markov chain's stepping class is expected to be approximately equal to the time from $t_\mathrm{start}$ (the start of the transition interval) to the transition. So, it should be on the same order but slightly shorter than the length of the transition interval. In comparison to this, the modulus of the second eigenvalue, $\lambda_\mathrm{dec}$, is associated to the \emph{mixing time} 
	\begin{align}
	    \MIX_F:=\frac{1}{1-\lambda_\mathrm{dec}}
	    \label{mixingtime}
	\end{align} within the stepping class. This mixing time is related to the time it takes to ``forget the initial condition'' while staying in the stepping class. We expect this time to be the time it takes to perform several steps (large scale oscillations in the data time series), so, on the order of several seconds.
	
	The magenta crosses in  \cref{metconclusion31} right
        highlight the boxes with center $X_{\mathrm{emp},i}$ for which
        $\MET_i$ are less than the overall mean
        $\MET_F=\frac{1}{1-\lambda_1}$. We observe that these magenta
        crosses are not uniformly  spread around the unit circle but
        are mostly concentrated in a range of phases. We note that the escape time $\MET_F$ and its
        distribution $\mathbf{s}_F$ is guaranteed to exist due to the  ergodicity of the stepping class $F_\mathrm{step}$:
        \begin{align}\label{def:sf}
          s_{F,i}=\lim_{N\to\infty}P(X_N=F_i\,|\,X_0=k\mbox{\ and\ }X_N\in F)\mbox{\ for all $k\in F$}
        \end{align}
        (where we use the notation $F_i$ to indicate the $i$th element
        of transition set $F$ in the Markov chain state space
        $\mathcal{J}_\mathrm{emp}$). In particular, the stationary
        distribution is reachable from all points in the transition set.

\subsection{Method summary}
\label{sec:method:sum}
In summary, for each part of the empirically measured force time series where experimenters (domain experts) have flagged a freezing event we proceed in the following way:
\begin{enumerate}
\item Select a transition interval $[t_\mathrm{start},t_\mathrm{end}]$ such that the data indicates stepping at $t_\mathrm{start}$ but freezing at $t_\mathrm{end}$ and such that several steps are included. 
\item Embed the scalar time profile $x(t)$ into the unit circle in the
  complex plane using the Hilbert Transform and scaling \eqref{scaling} to obtain $X(t)\in\mathbb{C}$ for $t\in[t_\mathrm{start},t_\mathrm{end}]$.
\item Subdivide the unit circle into $P\times Q$ boxes along polar coordinates ($P$ radial annuli of equal thickness, and $Q$ cones, see \cref{subdivide} for illustration).
\item Generate an empirical discrete-time discrete-space Markov chain
  transition matrix $A$ using \eqref{def:transitionmatrix}. This
  transition matrix is typically supported only on a subset of the
  $PQ$ boxes. Define the centers of these boxes $X_{\mathrm{emp},i}$,
  with angle coordinate $\psi_{\mathrm{emp},i}=\psi(X_{\mathrm{emp},i})$ for some
  $i\in\{1,\ldots,PQ\}$.
\item Identify the communicating classes for $A$ and test if one can
  split them into a transition set $F$ (including a stepping class,
  containing most of the states) and an absorbing set $E$.
\item Determine mean escape times $\MET_i$ from $F$ into $E$ when
  starting from any state $i$ in $F$. Define $X_\mathrm{min}$ as the coordinates of the state $i_\mathrm{min}\in F$ with minimal $\MET_i$ (the angle  coordinate is called $\psi_\mathrm{min}=\psi(X_{\min})$).
\end{enumerate}

\cref{discritise size} investigates how robust (as far as we tested for the data available) our results are with respect to the choice of $[t_\mathrm{start},t_\mathrm{end}]$ and discretization parameters $P$ and $Q$.

\section{Illustration of the methodology}
\label{sec:data}
\begin{figure}[ht]
	\centering
	\includegraphics[width=\textwidth]{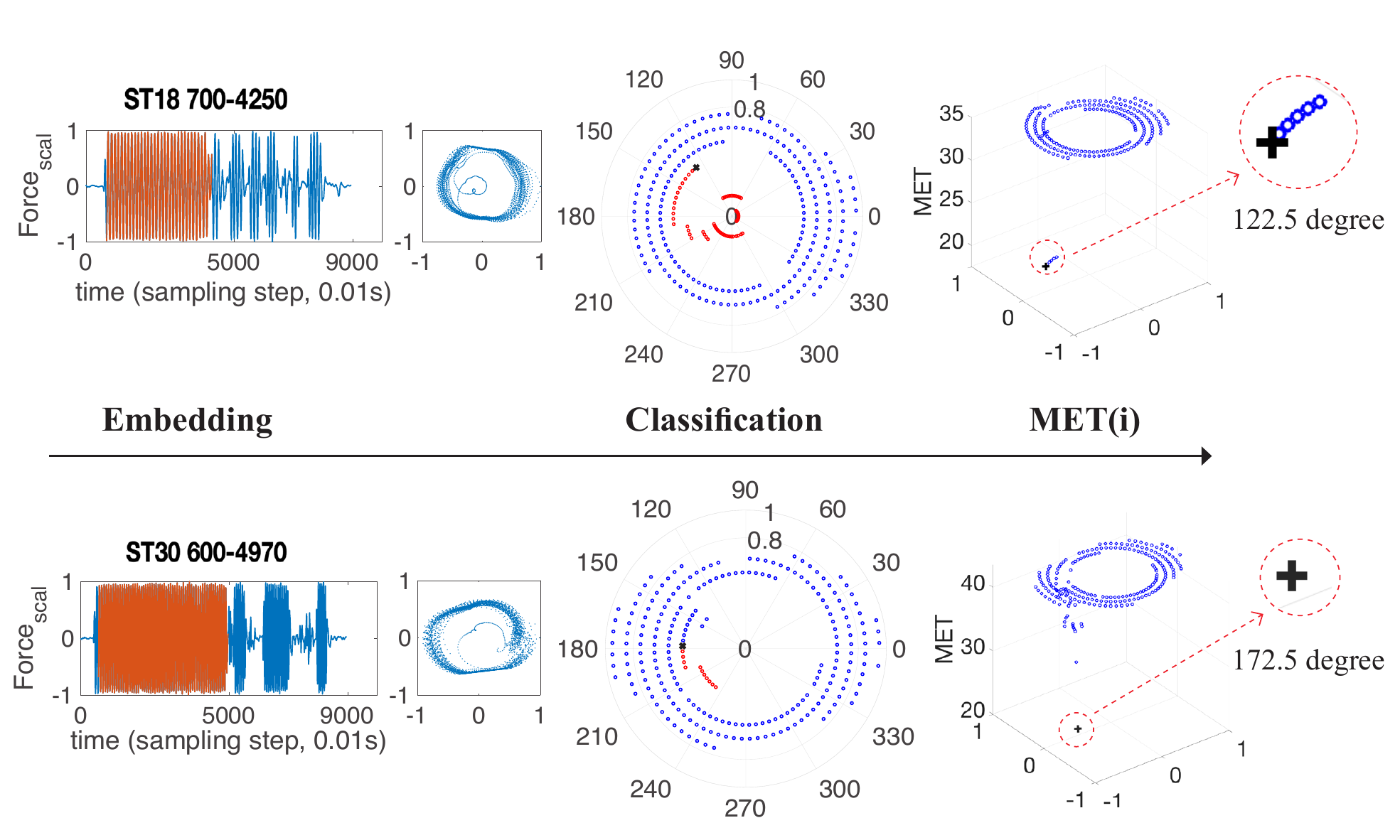}
	\caption{Illustration of the methodology for two subjects, ST18 and ST30. Step 1: choose the transition interval $[t_\mathrm{start}, t_\mathrm{end}]$  for embedding. Step 2: discretize complex plane into boxes along polar coordinate axes, obtain empirical transition probabilities and find communicating classes of resulting Markov chain. Step 3: determine $\MET_i$,  the mean escape time, from transition set $F$ for each box.}
	\label{comparison}
\end{figure}

To illustrate the methodology of using $\mathrm{MET}_i$ estimates to find preferred transition phases, we choose time series of two different patients. \cref{comparison} presents a graphical summary of the proposed methodology. 

The top panel starts with the time series from ST18, data set 1. In this case we use $[700,4250]$ as the transition interval, which amounts to a total of $35.50$\,s where the sampling time step is $0.01$\,s between sampling points. After Hilbert embedding and discretization, all the points can be classified into transition set $F$ or absorbing set $E$. There are a total $3455$ points in transition set $F$, which corresponds to a time interval of $34.55$\,s. Therefore the empirical transition time in the data for this transition interval is $34.55$\,s. The mean escape time $\MET_F$ for the resulting Markov chain is $34.4139$\,s (obtained via $\MET_F=1/(1-\lambda_1)$ following equation \eqref{metbylambda1}), so is roughly the same, as expected. The mixing time is $7.7675$\,s (obtained via $\MIX_F=1/(1-\lambda_\mathrm{dec})$, following equation \eqref{mixingtime}). The preferred phase given from $\psi_{\mathrm{min}}$ is $122.5^\circ$. The dominant eigenvalues for $A_F$ are $\lambda_1 = 0.9997$ and (in modulus) $\lambda_\mathrm{dec} = 0.9987$ after rounding to four decimal places.  

The bottom panel depicts time series from ST\,30, data set 1. In this case we use $[600,4970]$ as the transition interval, which amounts to a total of 43.70s with time step of $0.01$\,s between sampling points. There are a total 4352 points in transition set $F$, which correspond to $43.52$\,s, therefore for this transition interval the empirical transition time in the data is  $\MET_F=43.52$\,s.  The METs of the Markov chain is $43.1010$\,s (again, as expected, very close). The mixing time $\MIX_F$ is $5.8854$\,s. The preferred phase given from $\psi_{\mathrm{min}}$ is $172.5^\circ$. $\lambda_1 = 0.9998$ and $\lambda_\mathrm{dec} = 0.9983$ after rounding to four decimal places. 

\section{Transition phases for freezing events from stepping data}\label{sec:phase:prediction}

\begin{figure}[ht]
	\centering
	\includegraphics[width=0.9\textwidth]{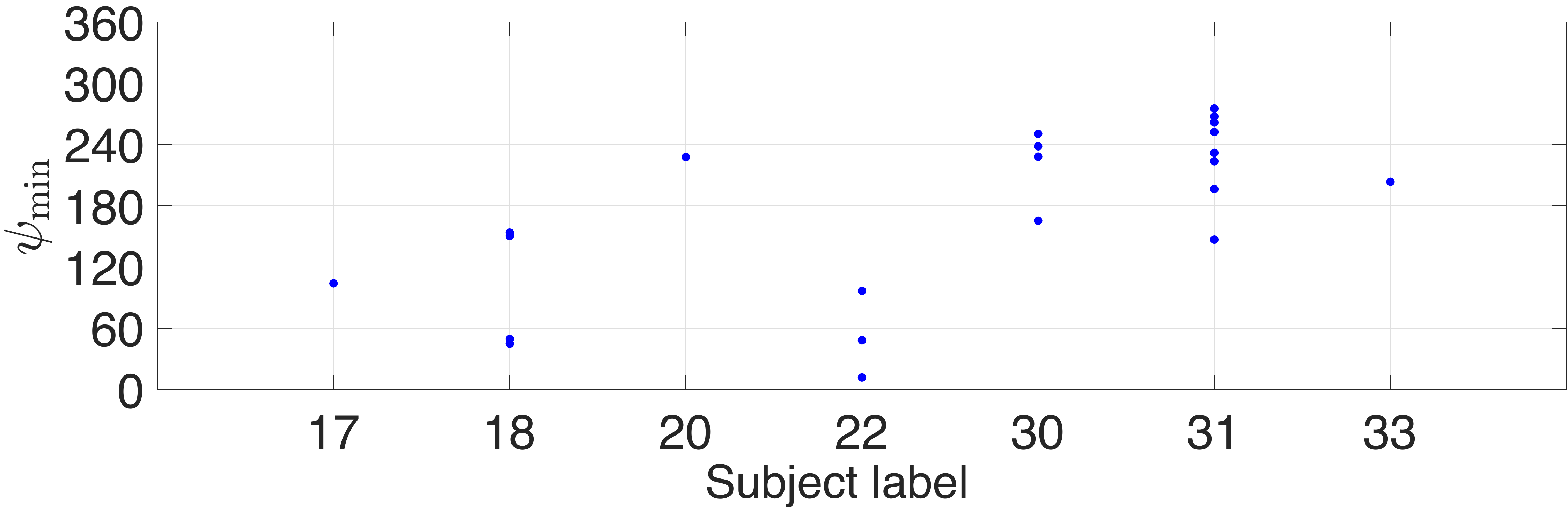}
	\caption{Phase $\psi_\mathrm{min}$ (preferred escape phase, in degree) of polar box $i$ with minimal expected escape time $\MET_i$ of all subjects showing distinct stepping and freezing episodes in the data set. Number of events: 1 (ST17), 4 (ST18), 1 (ST20), 3 (ST22), 4 (ST30), 8 (ST31), 1 (ST33)}
	\label{metforall}
\end{figure}
	\begin{figure}[ht]
		\centering
		\includegraphics[width=0.9\textwidth]{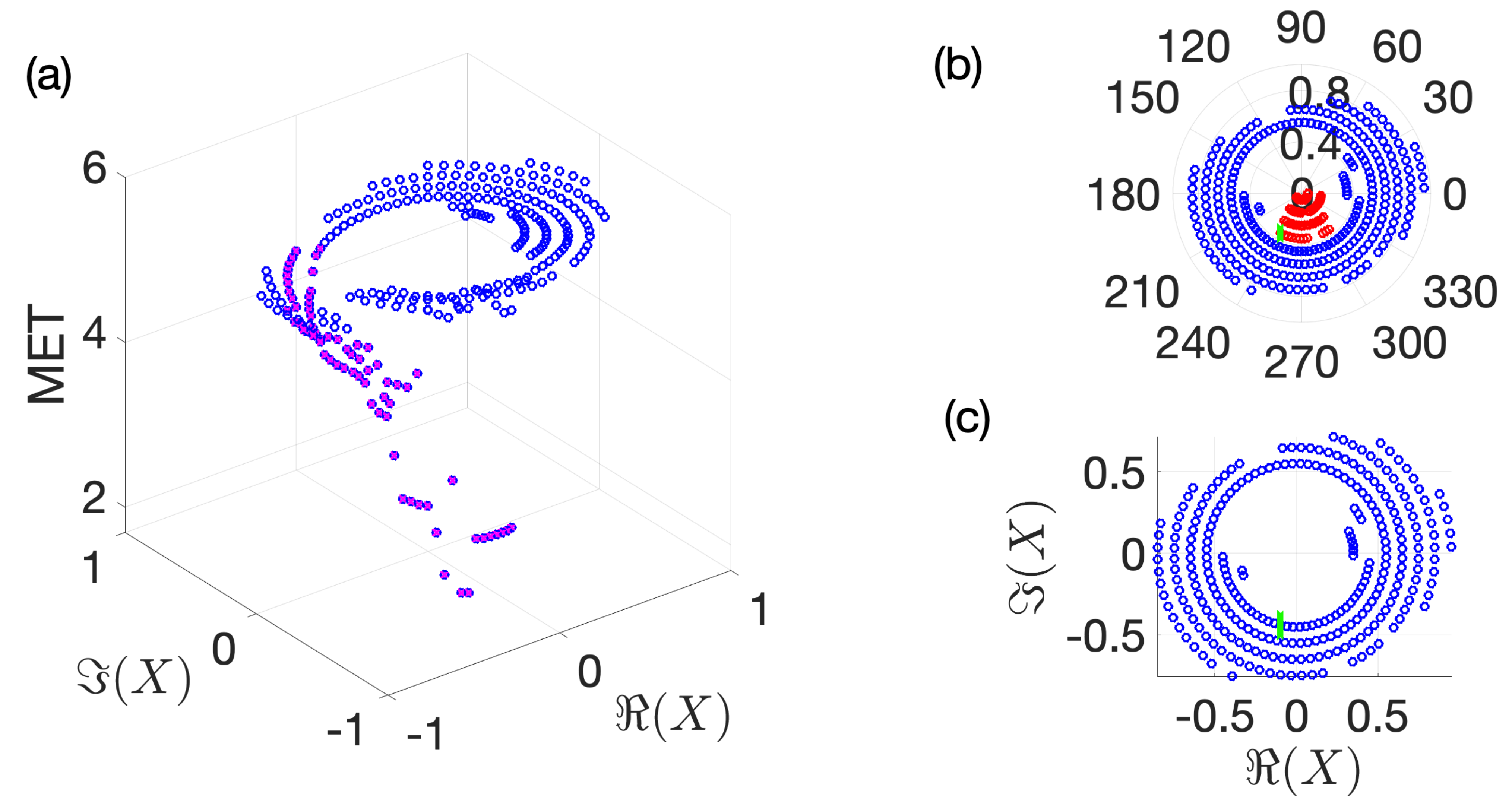}
		\caption{Communicating classes and mean escape times from transition set when combining all freezing episodes for subject ST31. Panel (a) shows the mean escape time calculated from each starting state on transition set $F$, the states marked in magenta crosses represent the MET value of these states are below than \eqref{metbylambda1}. The preferred angle is marked as black '+' in Panel (c). Panel (b) shows the state space of the trajectory after discretization and subdivision with $(p,q)=(0.1,5^\circ)$. The red dots represents the states in absorbing set $E$, the black cross represents the $\psi_\mathrm{tr}$ which is the first states that the Markov chain touches the absorbing set $E$. The blue dots represents the states in transition set $F$.}
		\label{combined_event31/1}
	\end{figure}
     \cref{metforall} shows the preferred transition phases $\psi_\mathrm{min}$ for all events present in the experimental data. The $x$-axis denotes the patient number. We observe that the values for $\psi_\mathrm{min}$ corresponding to different subjects appear to be different. 
    These observations suggest that although there might be a preferred transition phase for an individual subject it is not clear whether this holds true between subjects. Individual characteristics influencing the phase angle may be height, weight and the length of legs, etc.
    The limited number of repeating FE per subject in the data set does not allow us to be conclusive and reject the nullmodel. Nevertheless, \cref{metforall} indicates a non-uniform distribution of phases and hence the possibility of phase dependence of the dynamics underlying the transitions from stepping into freezing we have hypothesized in this paper. \cref{combined_event31/1} illustrates the proposed methodology applied to a combination of all events for patient \#31. We observe that the preferred phase indicated by \cref{combined_event31/1} is consistent with the phase angles calculated independently and presented in \cref{metforall}.

\section{Conclusions and Outlook}
\label{sec:conclusions}

In this paper we propose a methodology for timing the transition into freezing (of gait) and locating this transition in a reconstructed phase space. The methodology combines nonlinear time series analysis and mathematical modelling. We apply the developed methodology to real patient data collected as part of stepping in place experiments. Our approach is patient-specific and thus capable of studying properties of the FOG phenomenon on an individual level basis and potentially applicable on an individual level in the context of personalized interventions in the future \cite{lewis2022stepping}. Specifically, our analytical estimate $\frac{1}{1-\lambda_1}$ (equation \ref{metbylambda1}) for the time it takes a patient to transition into freezing given a segment of stepping time series data immediately preceding a FE could be used in future algorithms for online FE prediction that could be transformtive for the quality of life of patients with Parkinson's disease as it would enable early warning signals to be calibrated (for each individual) and implemented in wearable devices and/or pressure sensing shoe insoles.
Although our data set does not contain information on the step length, the use of insoles in \cite{pardoel2021grouping, shalin2021prediction} allows the force generated during forward walking to be measured. In this regard the temporal characteristics evaluated here should still apply to forward walking. While measurements from wearable devices may be noisier and hence less accurate than the stepping data from Nantel \emph{et al.}, in cases when the ratio between signal and noise is sufficiently large, phase space reconstruction based on embedding methods would be appropriate. The decomposability of the resulting state space of the Markov chain into transition set and absorbing set will be a good a-posteriori criterion whether the chosen embedding dimension is appropriate. This an interesting prospect for future work. 

Several hypothesis underlying the FOG phenomenon based on mechanistic studies have been proposed in the literature \cite{gao2020freezing}. These involve: a threshold mechanism \cite{plotnik2012freezing}; an inference mechanism \cite{lewis2009pathophysiological}; a cognitive mechanism \cite{vandenbossche2013freezing}; and a decoupling mechanism \cite{jacobs2009knee}.
In section 3 presenting a nullmodel, we hypothesise that the transitions associated with spontaneous involuntary freezing episodes that are not triggered by any apparent external stimuli (as was the case in the experiments producing the data set we used) can be described as noise-induced escape in a bistable oscillator setting. In this sense, our modelling approach is phenomenological and could be most closely associated with the threshold model \cite{plotnik2012freezing} mentioned above. This model does not account for mechanisms driving freezing but rather the dynamical (geometric) properties of the transition in to freezing. In fact, the analysis we have carried out and properties we have defined could be equally applied to both, patients and healthy individuals, should there be available data.

The quantity $\psi_\mathrm{min}$ extracted from the data by our method identifies for each freezing episode a unique threshold-independent time and stepping phase after which the subject is committed to freezing in this particular episode. The definition of such a characteristic quantity permits classification and clustering of freezing episodes and subjects according to in-person-between-episode and between-person mean and variation of $\psi_\mathrm{min}$. A non-uniform and possibly patient-specific distribution of preferred escape phases $\psi_\mathrm{min}$ is still compatible with the hypothesis of noise-induced escape from a limit cycle as proposed as an underlying mechanism in \cref{sec:nullmodel}. A non-uniform distribution would merely provide evidence against the rotation invariance present in nullmodel \eqref{eq:HNFc}. In this case one can explore whether this quantity $\psi_\mathrm{min}$ is correlated with other characteristics of freezing or the disease (e.g., severity or frequency of freezing events, progression of the disease or effectiveness of therapies). If it is the case that there may be a preferred phase at which an individual patient enters into freezing is another interesting result in terms of future model development that accounts for this transition. It could lead to the development of intervention strategies such as early warnings of FE based on patient-specific phase information.  
Therefore, a potentially very important direction for future work is to systematically determine preferred transition phases by analysing transitions into FE not just in stepping but also in the case of walking and freezing. In principle the method permits one to construct an empirical transition matrix from several time series by merging the box counts in \eqref{def:transitioncount}. Then we can not necessarily expect a single (unique) entrance state $i_\mathrm{tr}$. There will still be a unique $\psi_\mathrm{min}$, however, several local minima may make the transferred transition angle more uncertain.

\section*{Acknowledgments}
JS gratefully acknowledges support by EPSRC Fellowship EP/N023544/1 and EPSRC grant EP/V04687X/1. KTA gratefully acknowledges the financial support of the EPSRC via grant EP/T017856/1 and the support of the Technical University of Munich – Institute for Advanced Study, funded by the German Excellence Initiative. AW would like to acknowledge helpful comments and suggestions by Zhongkai Tao, University of California, Berkeley and Congping Lin, Huazhong University of Science and Technology.

\section*{Data availability}
Full data sets and processing scripts are
available at the following link \url{https://figshare.com/s/a14be7360925639736ba}.

\appendix
\section{Dependence on discretization parameters and length of transition intervals}\label{discritise size}
\begin{figure}[ht]
		\centering
 		\includegraphics[width=0.9\textwidth]{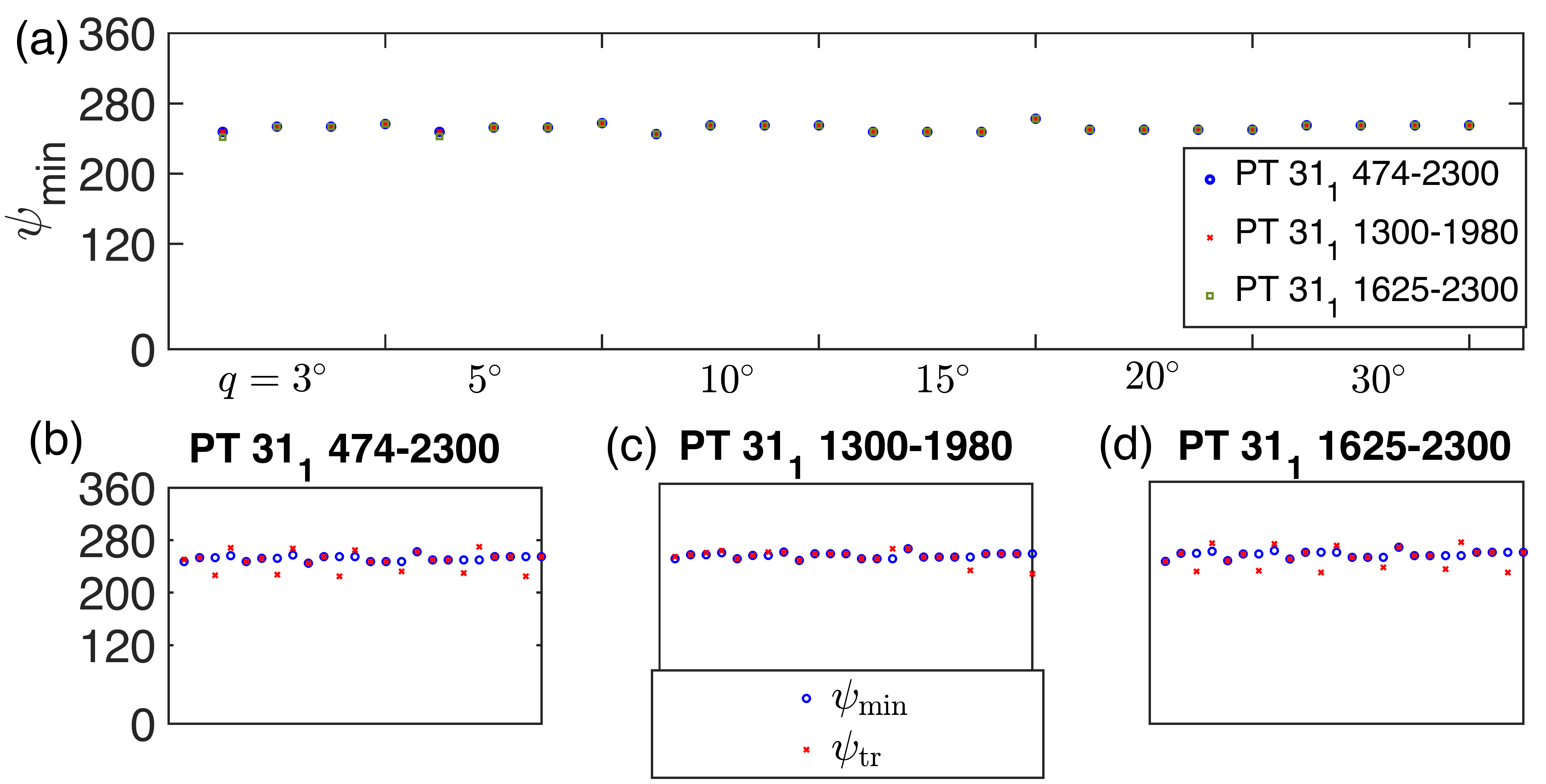}
		\caption{Panel (a): phase $\psi_\mathrm{min}$(in degree), where $\MET$ is minimal for different box sises and transition intervals.  The label on the $x$-axis shows box sizes $(p, q)$ for  $p\in\{0.05,0.1,0.15,0.2\},q\in\{3^\circ,5^\circ,10^\circ,15^\circ,20^\circ,30^\circ\}$. For each $p$, $q$ increases within its range. Blue circles, red crosses and green square markers correspond to transition interval lengths, as indicated in the legend. Panels (b)--(d) show dependence of $\psi_{\min}$ (blue circles) and $\psi_\mathrm{tr}$ (red crosses) on same box sizes (using same $x$-axis as panel (a)) and transition interval.}
		\label{metanalysis31}
	\end{figure}

\begin{figure}[ht]
\centering
	\includegraphics[width=0.8\textwidth]{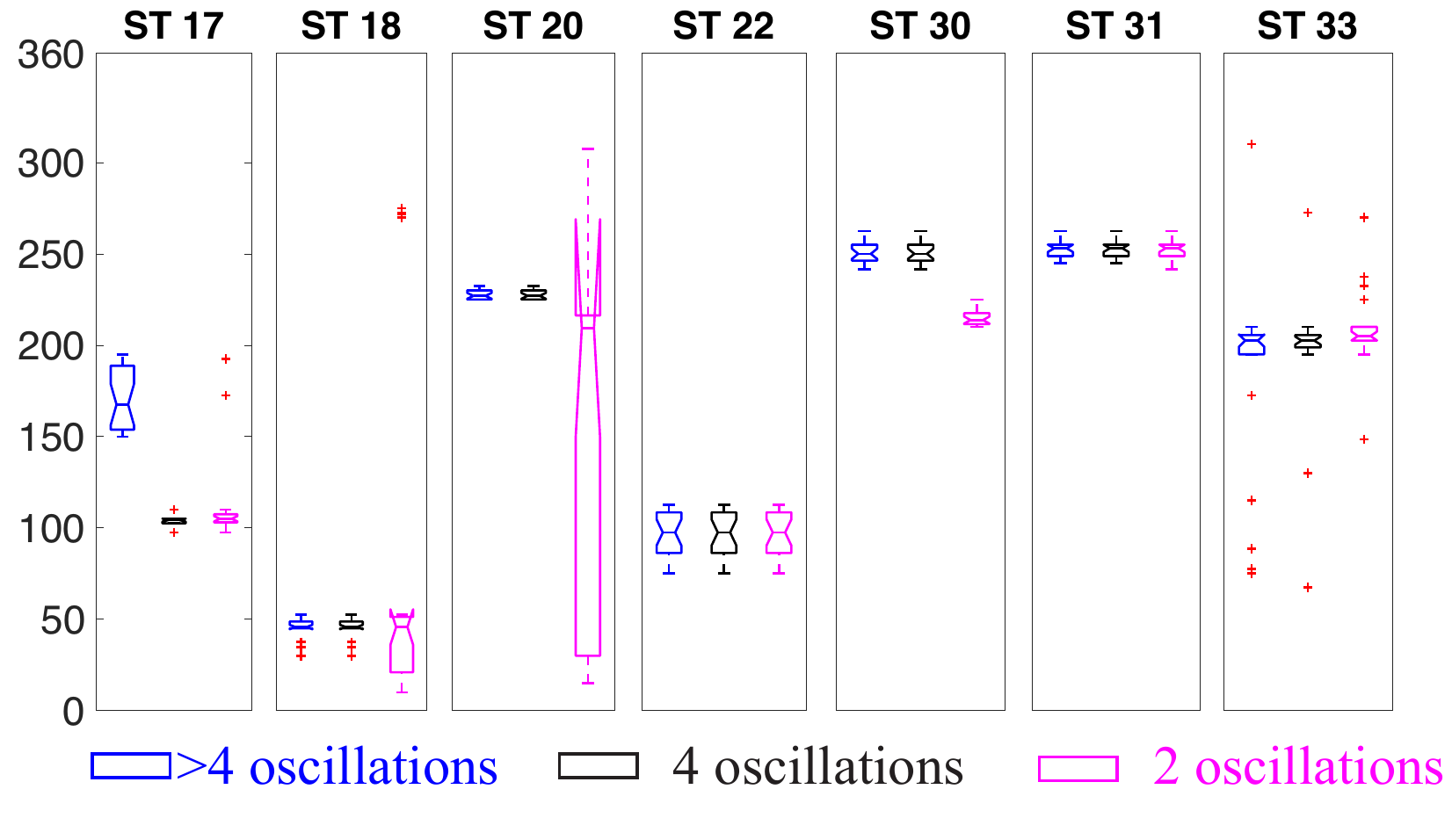}
\caption{Boxplot for phase $\psi_\mathrm{min}$ for same range of box sizes as in \cref{metanalysis31} and interval lengths for different subjects. The left boxplot (blue), middle boxplot (black) and right boxplot (magenta) of each patient represent 3 different choices of \emph{transition interval}: $[t_\mathrm{start}, t_\mathrm{end}]$ in the same time series that correspond to the maximal possible transition interval, approximately four oscillations and approximately two oscillations, respectively.}
\label{convar_analysis}
\end{figure}

\begin{figure}[ht]
		\centering
		\setlength{\abovecaptionskip}{0.cm}
		\setlength{\belowcaptionskip}{-0.cm}
 		\includegraphics[width=\textwidth]{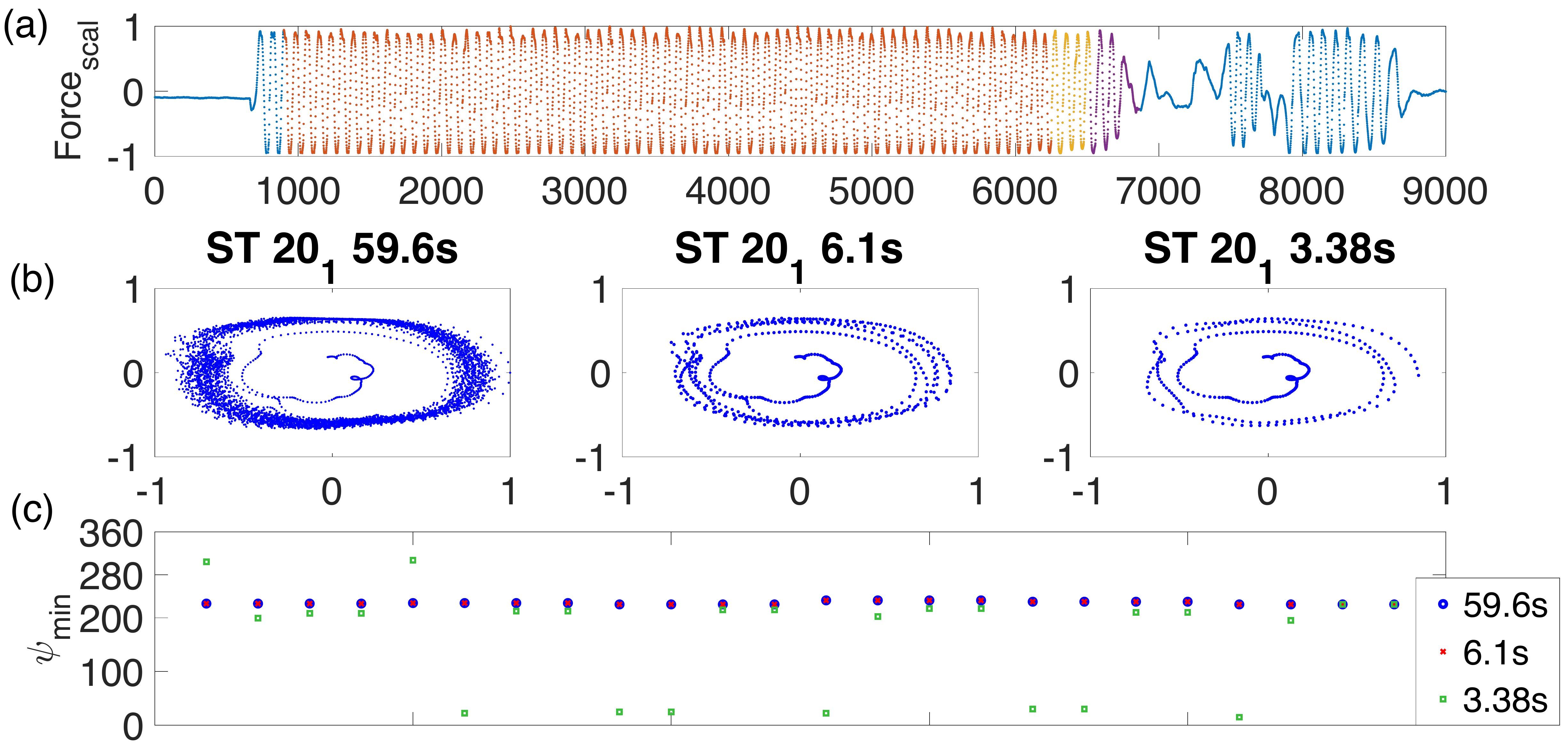}
		\caption{Panel (a) shows time series of freezing episode for subject ST20. Panel (b) shows the embedding trajectories for different transition interval lengths (left: many oscillations, middle: approximately $4$ oscillations, right: approximately $2$ oscillations) corresponding to the force data shown in Panel (a). Panel (c) shows the transition phase $\psi_\mathrm{min}$ for  different box
		sizes $p$ and $q$ (see \cref{metanalysis31}).}
		\label{ST20analysis}
	\end{figure}

 \Cref{metanalysis31,convar_analysis,ST20analysis} present a sensitivity analysis of the results with respect to various method parameters. 
Let us denote by $\psi_\mathrm{min}$ the phase (or angle) that corresponds to the box with minimal survival time $\MET_i$ in transition set $F$  (derived by \eqref{metbymatrix}).  \cref{metanalysis31}(a) shows  $\psi_\mathrm{min}$ for different discretization box sizes $q$, $p$, as introduced in \cref{sec:subdivison}. The $x$-axis in \cref{metanalysis31}(a) shows all combinations of $p$ and $q$ while the color and marker type encode different transition interval lengths (blue circles, red crosses and green squares, respectively). The transition intervals tested are: $[474,1980]$ (of length $15.06$\,s, containing more than four large-amplitude oscillations); $[1300,1980]$ (of length $6.8$\,s, containing approximately four large-amplitude oscillations); and $[1625,1980]$ (of length $3.55$\,s, containing approximately two oscillations) respectively. 
Panels (b--d) in \cref{metanalysis31} compare for the same range of box sizes  $(p,q)$ and transition intervals the sensitivity of $\psi_{\min}$ (blue dots, the phase determined by minimal $\MET_i$ in the transition set $F$) to the sensitvity of $\psi_\mathrm{tr}$ (red crosses, the phase determined by the first state reachable in the absorbing set $E$). We observe that the  phase $\psi_\mathrm{min}$ is less sensitive to box sizes and transition interval length than $\psi_\mathrm{tr}$.
\cref{convar_analysis} summarizes the sensitivity for further freezing episodes and different subjects as box plots for the phase $\psi_\mathrm{min}$, varying discretization box sizes over the same range as \cref{metanalysis31} within each box plot. The left, middle and the right box plot of each subject are for $3$ different transition interval lengths $[t_\mathrm{start}, t_\mathrm{end}]$: the transition interval for the left box plot contains more than $4$ stepping periods (large-amplitude oscillations), the transition interval for the middle box plot contains approximately $4$ oscillations, and the transition interval for the right box plot contains $2$ oscillations, respectively. The distribution of extracted $\psi_\mathrm{min}$ supports our general observation that, if the stepping class contains sufficiently many oscillations (approximately $4$), the resulting transition phase $\psi_{\min}$ is largely insensitive to our method parameters. 

\Cref{ST20analysis} shows in more detail the transition phases $\psi_\mathrm{min}$ for subject ST20 (\cref{convar_analysis}). In panel (b), in the case of 2 oscillations, the data in the stepping class is not sufficient to divide the state space into classes (transition set $F$ and absorbing set $E$), as they are too sensitive to the transition interval length in this example.  
The above sensitivity analysis indicates that $\approx4$ oscillations preceding a freezing episode in the transition interval are recommended in order to obtain robust phase estimate in this data set.

\section{Embedded time series for different subjects}
\label{sec:more:events}
\begin{figure}[ht]
		\centering
		\setlength{\abovecaptionskip}{0.cm}
		\setlength{\belowcaptionskip}{-0.cm}
		\includegraphics[width=\textwidth]{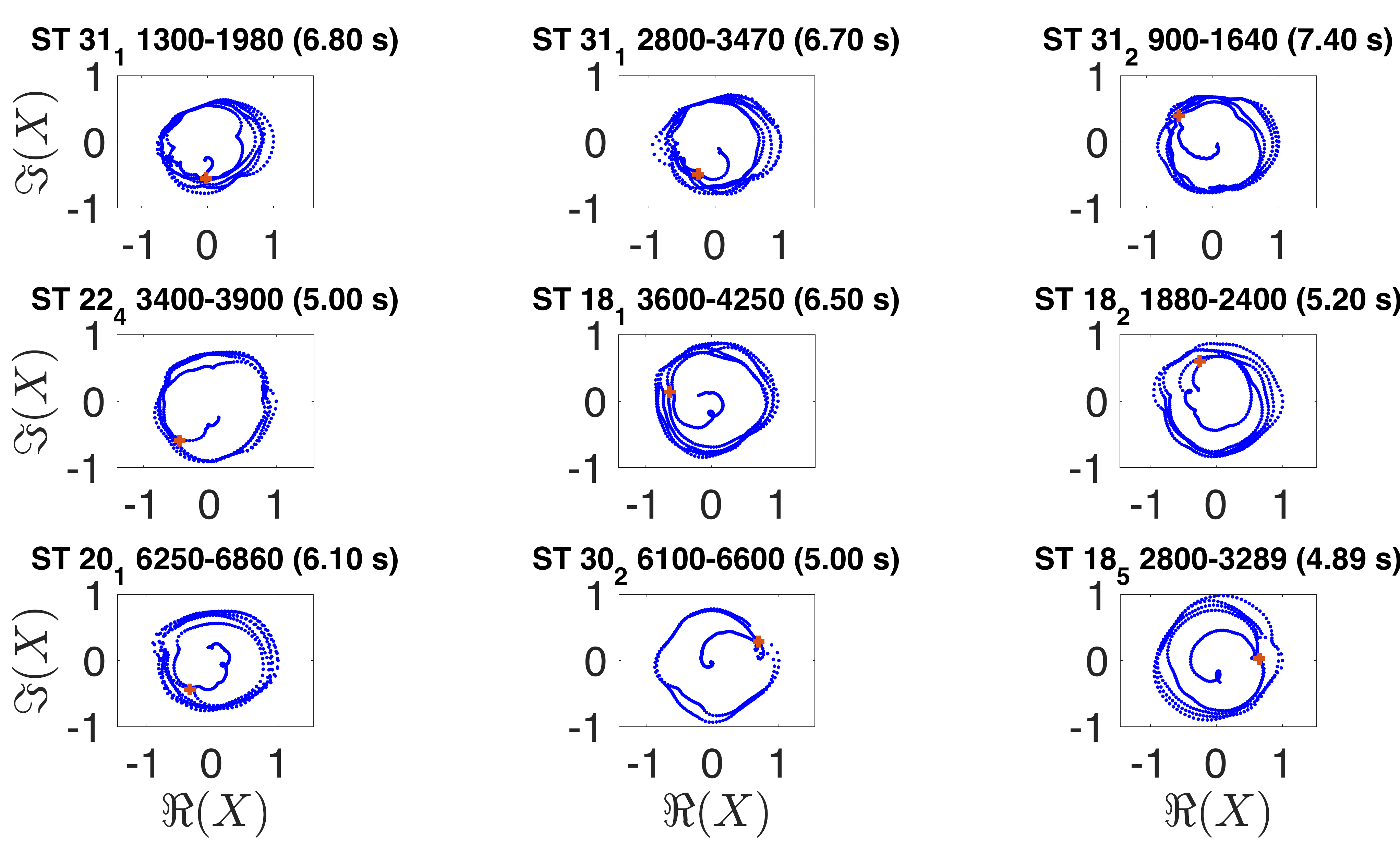}
		\caption{Embedded trajectories in $\mathbb{C}$ after scaling and Hilbert transform for $9$ selected events and transition intervals. Each panel title specifies subject number, subject trial data set and sampling step numbers and resulting length of transition interval in seconds. See also step~1 of \cref{comparison}.  Red crosses indicate the cartesian coordinates $X_\mathrm{min}$ of boxes with minimal escape time from transition set.
		\label{hilbert_polar1}}
	\end{figure}

\begin{figure}[ht]
		\centering
		\setlength{\abovecaptionskip}{0.cm}
		\setlength{\belowcaptionskip}{-0.cm}
		\includegraphics[width=\textwidth]{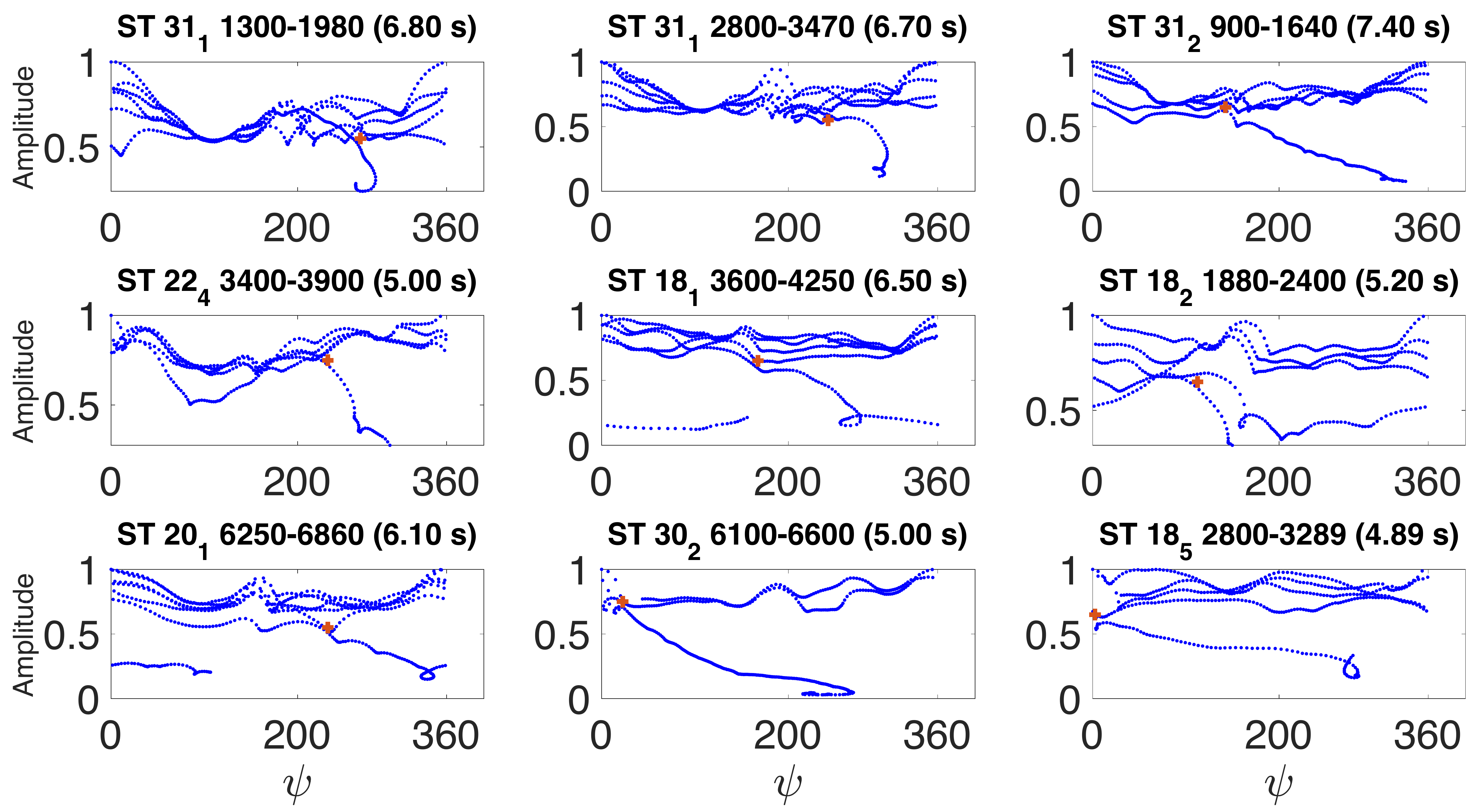}
		\caption{Same embedded trajectories as \cref{hilbert_polar1} in polar coordinates $(\psi,R)$, where $\psi$ is in degree and $R$ represents the amplitude.  Red crosses indicate the polar coordinates $(R_\mathrm{min},\psi_\mathrm{min})$ of boxes with minimal escape time from transition set.}
		\label{amplitude_degree1}
	\end{figure}

\cref{hilbert_polar1} shows the embedding of time series for different subjects with freezing events, obtained by applying the Hilbert Transform. 
\cref{hilbert_polar1} shows the same embedded trajectories as \cref{amplitude_degree1} in polar coordinates $(\psi,R)$, where $\psi$ is in degree and $R$ is scaled to $[0,1]$ after Hilbert Transform. 
	
\end{document}